\documentclass[reqno]{amsart}
\usepackage[utf8]{inputenc}

\usepackage[utf8]{inputenc}

\usepackage{sky}
\usepackage{graphicx}
\usepackage{color}
\usepackage{caption}
\usepackage{subcaption}

\newcommand{\area}{\mrm{area}}

\newcommand{\tail}{\tau}
\newcommand{\modarea}{\underbar{area}}

\title{Expanded regimes of area law for lattice Yang--Mills theories}

\author{Sky Cao}
\address{Department of Mathematics, Massachusetts Institute of Technology, Cambridge, MA 02139}
\email{skycao@mit.edu}

\author{Ron Nissim}
\email{rnissim@mit.edu}

\author{Scott Sheffield}
\email{sheffield@math.mit.edu}

\date{}

\begin{document}

\begin{abstract}
We extend the parameter regimes for which area law is proven for pure $\mathrm{U}(N)$ lattice Yang--Mills theories, in particular when $N$ is large.  
This improves on a classical result of Osterwalder-Seiler from 1978. To do so, we view the master loop equation as a linear inhomogeneous equation for Wilson string expectations, and then prove an a priori bound for solutions to the equation. The main novelty is in how we deal with the merger term in the master loop equation. This is done by introducing a truncated model for which the merger term is unproblematic, and then showing that the truncated model well approximates the original model.
\end{abstract}

\maketitle

\tableofcontents

\section{Introduction}\label{section:introduction}

The problem of quark confinement in four dimensions is one of the major open problems in Yang--Mills theory. For instance, it is mentioned as one of the natural extensions to the Yang--Mills millennium problem in \cite{Jaffe2006a}. The problem amounts to showing a rapid decay property for a family of observables known as Wilson loop expectations. This rapid decay property is known as {\em area law}, and was introduced by Wilson \cite{Wilson1974} as a potential explanation for the phenomenon of quark confinement (that is, single quarks are never observed in nature). The conjecture is that in four dimensions, area law should hold for non-Abelian Yang--Mills theories at all values of the inverse coupling constant $\beta$, which is a parameter that governs the amount of order/disorder in the system. Due to physical reasons, the main theory of interest is $\mrm{SU}(3)$ Yang--Mills theory, but it would already be mathematically interesting to show area law at all $\beta$ for any non-Abelian theory.

As far as we are aware, the only {\em prior} result showing area law in dimension four in some regime of $\beta$ is the classical result of Osterwalder-Seiler \cite{osterwalderseiler1978}, which works in the strong coupling regime (where $\beta$ is small). We briefly survey other relevant works. Fr\"{o}hlich-Spencer \cite{frohlich1982massless} showed that area law does not hold for the $\unitary(1)$ Yang--Mills theory in four dimensions at large enough $\beta$ (so that the non-Abelian assumption is necessary). In dimension three, area law is known to hold for the $\unitary(1)$ theory at all $\beta$ \cite{gopfertmack1981}, and then by a general result of \cite{FROHLICH1979}, this implies that the same is true for $\unitary(N)$ Yang--Mills theory for any $N \geq 1$. We remark that \cite{FROHLICH1979} also applies in general dimensions, but by tracking the constants in that paper, one sees that in general dimensions, it reproves the Osterwalder-Seiler result \cite{osterwalderseiler1978}. The works \cite{Chatterjee2019a, jafarov2016} showed area law in the large-$N$ limit, under certain conditions on $\beta$. The works \cite{durhuus1980connection, chatterjee2021probabilistic} gave sufficient conditions for area law to hold, in terms of strong enough mass gap assumptions.  

\begin{remark}
Although \cite{osterwalderseiler1978} is the strongest {\em prior} result, we should note that since this paper first appeared on the arXiv, there has been a {\em posterior} result: namely, the current authors in \cite{CNS2025b} were able to establish area law for a larger range of $\beta$ values (albeit with weaker constants for the exponential decay rates) using a completely different approach---one that combines dynamical techniques recently introduced in \cite{shen2023stochastic} with more classical {\em area-law-proof} techniques from \cite{durhuus1980connection}. Nonetheless, is is worth noting that several features distinguish the current paper from \cite{CNS2025b}. First, the proof in this paper is inspired by the ideas of {\em gauge string duality} (string trajectories, surfaces) and it has long been of independent interest to explore what results could be obtained from this duality.
Second, as noted above, the ``string tension'' bound obtained in this paper is better than the one in \cite{CNS2025b}, see item~\ref{item:dynamical-approach-comparison} in Remark~\ref{remark:after-thm}, so at least for now there is also some technical advantage to the string theoretic point of view. Finally, as further discussed below, it remains a major open problem to show area law for large $\beta$ and $d=4$ for arbitrary $N>1$ (e.g.\ $N=2$ or $N=3$) and to {\em control the corresponding string tension constants} in the large $\beta$ limit. Area law is known {\em not} to hold for large $\beta$ and $d=4$ when $N=1$, so any proof will have to use $N>1$ in a crucial way. The dynamical and string-theoretic approaches provide different answers to the question of {\em why} one might expect taking $N>1$ to make a difference: in the dynamical approach of \cite{shen2023stochastic} and \cite{CNS2025b}, the difference is the positive curvature of groups like $\mathrm{SU}(N)$ when $N>1$ and the implications for the mixing rates of a type of Glauber dynamics on the configuration space. In the surface sum or string trajectory formulation, the difference is that when $N>1$, there is a bias toward surface configurations with {\em many small and low genus components}. Introducing a large Wilson loop effectively forces the existence of a large ``defect'' (e.g.\ a large surface or string trajectory) and one expects to pay a ``cost'' proportional to the size of the defect (although controlling this cost becomes more challenging when $\beta$ is large and $N$ is small, see Remark \ref{remark:future-directions}).
\end{remark}


In recent years, there have also been a number of works in probability studying other aspects of Yang--Mills, see e.g. \cite{shen2018,  GS2023, shen2022new, shen2023stochastic, chevyrev2023invariant, CPS2023, CG2025, AN2023, Cha2024, BCSK2024, CC2024, SSZ2024, SZZ2024, BCSK2025} for results on the lattice, and \cite{che19, CCHS2D, cao2021yang, cao2021state, CCHS3D, ChevyrevSurvey, BC2023, BC2024, CS2025} for results in the continuum. For additional discussion and even more references about area law and Yang--Mills more broadly, see the survey \cite{chatterjee2016a}.

\begin{figure}
\centering
\begin{subfigure}{.5\textwidth}
  \centering
\includegraphics[width=.8\linewidth, page=2]{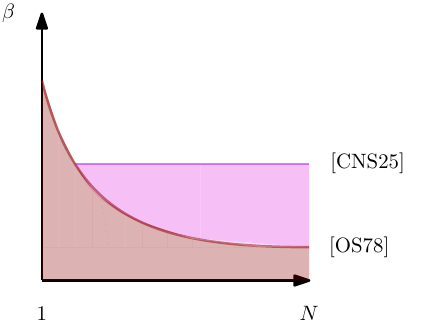}
\end{subfigure}%
\begin{subfigure}{.5\textwidth}
  \centering
\includegraphics[width=.8\linewidth, page=3]{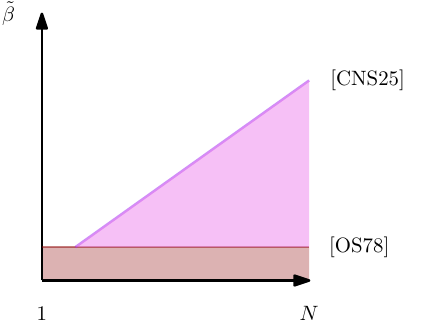}
\end{subfigure}
\caption{{\bf Left:} The two curves correspond to the $\beta$ thresholds below which area law is proven to hold. The dark red curve corresponds to the Osterwalder-Seiler result \cite{osterwalderseiler1978}, which shows that taking $\beta \ll_d 1/N$ (i.e., $\beta \leq c(d) /N$ for some small dimensional constant $c(d)$) suffices. We show that taking $\beta \ll_d 1$ suffices, which corresponds to the violet line. In particular, our threshold does not decay with $N$. Due to various choices of dimensional constants in our argument, our threshold may initially start out lower than the Osterwalder-Seiler threshold, but once $N$ is large enough, our threshold is better.
{\bf Right:} If one uses the {\em normalized trace} in defining lattice Yang--Mills (i.e.\ if one uses {\em eigenvalue mean} instead of {\em eigenvalue sum}; see Section \ref{section:preliminaries}), then the parameter $\beta$ is effectively replaced by $\tilde{\beta} = N \beta$. With this scaling, \cite{osterwalderseiler1978} proves area law for $\tilde{\beta} \ll_d 1$ (dark red line), while we establish area law for $\tilde{\beta} \ll_d N$ (violet line). This scaling is natural for some purposes, see Remark \ref{remark:trace-vs-normalized-trace}. From either perspective, it is an open problem to {\em further} extend the purple region upward and leftward. Treating small $N$ and large $\beta$ (or equivalently, large $\tilde \beta$) is of particular interest.}
\label{figure:os-comparison}
\end{figure}

In this paper, we improve on the classical Osterwalder-Seiler result \cite{osterwalderseiler1978} by showing area law for pure $\unitary(N)$ Yang--Mills theories in an extended regime, in particular when $N$ is large. For a visualization of the difference between the two regimes, see Figure \ref{figure:os-comparison}. Our main theorem is as follows. Any unfamiliar notation will be defined in Section \ref{section:preliminaries}.

\begin{theorem}[Area law]\label{thm:area-law-main}
Consider lattice Yang--Mills theory with gauge group $G = \unitary(N)$. For all $d \geq 2$, there exists $\beta_0(d) > 0$ such that for all $\beta \leq \beta_0(d)$ and all $N \geq 1$, for any rectangular loop $\ell$ belonging to a finite lattice $\Lambda$, we have that 
\begin{equs}
|\langle W_\ell \rangle_{\Lambda, \beta, N}| \leq C_{1, d} C_N^{|\ell|} \exp(-C_{2, d} \area(\ell)),
\end{equs}
where $|\ell|$ denotes the perimeter of $\ell$, and $C_{1, d}, C_{2, d}$ are constants depending on $d$, and $C_N$ is a constant depending on $N$. 
\end{theorem}

\begin{remark}\label{remark:after-thm}
We make the following remarks regarding Theorem \ref{thm:area-law-main}.
\begin{enumerate}
\item The result is really a result about large enough $N$, because if $N$ is bounded, then we could apply the Osterwalder-Seiler result \cite{osterwalderseiler1978}. Thus throughout this paper, one should think of $N$ as larger than some dimensional constant. One main point behind Theorem \ref{thm:area-law-main} is that when $N$ becomes large, it should become easier to prove area law. This is due to various reasons which we will discuss later. Our proof manages to take advantage of these reasons, whereas the Osterwalder-Seiler result does not, which is why it gives a threshold which becomes worse with $N$.
\item The result is uniform in the lattice $\Lambda$, and thus the same estimate will hold for any infinite-volume subsequential limit. In fact, in our current regime, one would expect the infinite volume limit to be unique, as was shown for $\mrm{SO}(N), \mrm{SU}(N)$ in \cite{shen2023stochastic}. One would expect that their proof can be extended to $\mrm{U}(N)$, althought currently this is not written anywhere, to our knowledge.

\item As long as $\ell$ is a (sequence of) rectangular loop(s) with both side lengths tending to infinity, the factor $C_N^{|\ell|}$ will be swamped by the factor $\exp(-C_{2, d} \area(\ell))$, because in this case $|\ell| \ll \area(\ell)$. Thus for such sequences of rectangular loops (which indeed are of primary interest), the area decay becomes more manifest. 

\item\label{item:dynamical-approach-comparison} (Comparison with the dynamical approach) The area law upper bound proved in \cite[Theorem 1.6]{CNS2025b} at present exhibits poor dependence on the parameter $N$. That is, the analog of the constant $C_{2, d}$ (which governs the rate of area law decay) from Theorem \ref{thm:area-law-main} actually decays with $N$ in \cite[Theorem 1.6]{CNS2025b} (see the footnote on \cite[Page 6]{CNS2025b}, which is related to \cite[Remark 4.12]{shen2023stochastic}). On the other hand, the constant $C_{2, d}$ in Theorem \ref{thm:area-law-main} does not decay with $N$. Thus, in the regime $\beta \leq \beta_0(d)$ of the current paper, the upper bound of Theorem \ref{thm:area-law-main} is better than the corresponding upper bound of \cite[Theorem 1.6]{CNS2025b}.

\item We expect that our proof approach should work for various other groups such as $\mrm{SU}(N)$, and perhaps even $\mrm{SO}(N)$, after various minor adjustments. In this paper, we focus on the most concrete case of $\mrm{U}(N)$, where the master loop equation (to be discussed soon) is the simplest.
\end{enumerate}
\end{remark}

In order to prove Theorem \ref{thm:area-law-main}, we introduce new techniques for analyzing lattice Yang--Mills theory. In particular, we develop a way of analyzing the finite-$N$ Master loop equation, which is a recursive equation satisfied by Wilson string expectations. In the past, this equation has only been understood in the large $N$-limit, see \cite{Chatterjee2019a, jafarov2016, chatterjee2016}. In this limit, the solution to the equation may be expressed as an absolutely convergent expansion (essentially by Picard iteration), while at finite $N$, the same does not hold, which forms an immediate barrier to understanding the finite-$N$ equation.

To circumvent this issue, we compare the lattice Yang-Mills theory to a ``truncated model'' for which solutions to the ``truncated Master loop equation'' may indeed be expressed as a convergent expansion. We then extend this argument back to the original model, by essentially showing that the truncated model very well approximates the original model. We sketch the argument in more detail in the paragraphs following the next two remarks.

As discussed more thoroughly in Remark \ref{remark:surface-exploration}, one can interpret our approach as implementing a ``surface exploration process'' for the surfaces constructed in \cite{CPS2023}. In particular, our argument can be thought of as analyzing lattice Yang--Mills theory via a geometric representation for its correlation functions (i.e.\ Wilson string expectations). For more discussion on this general approach towards understanding statistical mechanical models, see \cite[Section 1.1]{BCSK2025}. We hope that our approach can be further developed, for instance to handle the future directions discussed in Remark \ref{remark:future-directions}.

\begin{remark}[Future directions]\label{remark:future-directions}
It would be natural to try to extend Theorem \ref{thm:area-law-main} to the case where $\beta$ is allowed to depend linearly in $N$, e.g. $\beta \sim c_d N$, where $c_d$ is a small dimensional constant. In other words, the natural next step is to further develop the approach of the present paper to reprove the main result of \cite{CNS2025b}. This is the regime of the large-$N$ limit, in that upon sending $N$ to infinity, we obtain the limiting theory described in \cite{Chatterjee2019a, jafarov2016, BCSK2024}. While much is known about the limiting $N \toinf$ theory, comparatively less is known about more quantitative questions related to rates of convergence to the limit. We remark that while the result \cite{CNS2025b} covers the regime $\beta \sim c_d N$, at present it does not imply such quantitative statements, due the poor dependence on $N$ of certain constants, which was discussed in item~\ref{item:dynamical-approach-comparison} in Remark~\ref{remark:after-thm} (in principle it may be possible to refine the arguments of \cite{shen2023stochastic, CNS2025b} to obtain such quantitative statements, but so far this has not yet been done). 


If we are able to extend the approach of the current paper to the linear regime $\beta \sim c_d N$, then we expect that by necessity, we will be able to obtain quantitative estimates on rates of convergence. In order to obtain the extension, it appears that new ideas are needed even for the truncated model. If one wants to still perform some type of contraction mapping argument, then the main question is in finding an ansatz for the norm appearing in the contraction estimate \eqref{eq:contraction}. Presumably, it needs to account for certain cancellations that occur in the regime $\beta \sim c_d N$, which currently we are not utilizing. To put it differently, if a proof were formulated in terms of surface sums, it might require showing that a certain sum of surface weights is small, even though the sum of the absolute values of those weights is large; to this end, one would have to show that the positive and negative weights very nearly cancel each other out. This issue may become even more delicate when one goes beyond the $\beta \sim c_d N$ regime, which would be the natural next step after achieving the linear regime $\beta \sim c_d N$.
\end{remark}

Next, we discuss some of the main ideas in the proof of Theorem \ref{thm:area-law-main} in more detail. The main property we use is that the Wilson string expectations satisfy a linear inhomogeneous equation, known as the master loop equation. This equation appeared in various forms in \cite{Chatterjee2019a, chatterjee2016, jafarov2016, shen2022new, CPS2023, AN2023, BCSK2024}, as well as in many physics references (see \cite{Chatterjee2019a}). The approach is then to prove an a priori estimate for solutions to this equation, which amounts to proving a contraction estimate for the linear operator appearing in the equation. More explicitly, we are just using the simple fact that if we have a matrix $G$ and vectors $f_*, g$ such that
\begin{equs}
f_* = G f_* + g,
\end{equs}
and moreover $G$ satisfies the contraction estimate
\begin{equs}\label{eq:contraction}
\|G f\| \leq \frac{1}{2} \|f\| \text{ for all vectors $f$},
\end{equs}
where $\|\cdot\|$ is some appropriate norm, then we automatically have that (assuming of course that $\|f_*\|, \|g\| < \infty$)
\begin{equs}
\|f_*\| \leq 2 \|g\|.
\end{equs}
In the context of Yang--Mills, one should think of $f_*$ as a function mapping strings to Wilson string expectations, and the norm $\|\cdot\|$ is defined so that if $\|f_*\| \lesssim 1$, then area law decay holds.

This PDE viewpoint was first introduced in the context of Yang--Mills by Chatterjee in \cite{chatterjee2016}, and subsequently applied in \cite{chatterjee2016, jafarov2016}. All these cited works dealt with the large-$N$ limit (or expansions around the large-$N$ limit), because in this case the master loop equation simplifies, and so it is more manageable to prove the contraction estimate \eqref{eq:contraction}, which then directly leads to area law.

When $N$ is finite, the master loop equation becomes more complicated, due to the presence of the so-called ``merger term". These merger terms come with a prefactor of $N^{-2}$, which is why they are not present in the large-$N$ limit, and also why we expect that area law is easier to prove when $N$ becomes large.

Nevertheless, even with the $N^{-2}$ prefactor, this merger term forms an essential obstruction towards proving the contraction estimate \eqref{eq:contraction}. In particular, the previous arguments in \cite{Chatterjee2019a, chatterjee2016, jafarov2016, BCSK2024} do not directly apply. One of the main novelties in the present work is a new approach towards managing this problematic term, at least when $N$ is large. To do so, we first introduce a truncated version of lattice Yang--Mills, which is defined precisely so that the merger term becomes unproblematic. Thus the contraction estimate \eqref{eq:contraction} can be proven for this truncated model, from which area law follows. On the other hand, in our regime of parameters, one expects the truncated model to approximate the original model very well. Therefore, one hopes to be able to extend some form of the contraction estimate to the original model. The main idea in the extension step is that it should be extremely rare for the merger term to be problematic, and thus instead of trying to prove the contraction estimate in all situations, we content ourselves with proving the contraction estimate in the ``typical" situation. In the atypical scenarios where we come across a problematic merger term, we a priori know that this is quite rare, and so such a term already comes with a huge cost.

In slightly more detail, the truncated model is defined by imposing an upper bound on the number of copies of each plaquette in the system. As a consequence of this upper bound, the number of possible merger terms itself becomes bounded, which ultimately is why these terms become unproblematic. On the other hand, the imposed upper bound is still large enough so that it is extremely rare to see plaquettes with more copies than the upper bound, which in the end is the reason why we are able to obtain the result for the original model. In this discussion, we have neglected one key idea which is needed to handle this extension step, because the discussion of this idea would perhaps be too detailed in current introductory section. Thus, we include an even more detailed discussion of the extension step near the beginning of Section \ref{section:original-model}.

At a high level, this strategy for extending the contraction estimate to the original model is analogous to following way of refining the second moment method in probability. In various problems of probability, the main thing to show is that for some given random variable, the second moment is of the same order as the square of the first moment. However, there are certain situations where this is not the case, because the second moment is much larger. On the other hand, if one can show that this largeness arises from events of very low probability, then one can hope to make a modified version of the second moment method work, where one computes the second moment after excluding these low probability problematic events. 

\begin{remark}[Comparison to the proof of \cite{osterwalderseiler1978}]
To compare with the proof strategy of Osterwalder-Seiler \cite{osterwalderseiler1978}, there the authors set the truncation level to be at one plaquette, in that they assumed that $\beta$ was small enough so that it is very rare to see even a single plaquette. The corresponding truncated model is then a system with no plaquettes, for which area law trivially holds.

Our refinement of \cite{osterwalderseiler1978} arises because we allow our truncation level to depend on $N$; as $N$ increases, so does our truncation level. The additional input needed is an area law result for more general truncated models, where the truncation level is allowed to be much larger than $0$. As previously discussed, we are able to show such a result via a contraction argument, recall \eqref{eq:contraction}.
\end{remark}

\begin{remark}[Relation to the Dyson-Schwinger equation]
The master loop equation can also be thought of as the Yang--Mills analog of the Dyson-Schwinger equation for matrix models. The Dyson-Schwinger equation has been an important tool for studying matrix models and interacting particle systems, see e.g. \cite{borot2013asymptotic,guionnet2015asymptotics, BGK15, guionnet2019asymptotics,GH24}. However, similar to the previous works \cite{Chatterjee2019a, chatterjee2016, jafarov2016, BCSK2024} in the lattice Yang-Mills setting, in these other settings the Dyson-Schwinger equation was used to recover $1/N$ asymptotic expansions or other asymptotic properties for $N \to \infty$. In contrast, we are able to prove a non-asymptotic result, i.e. area law at finite $N$.
\end{remark}

\begin{remark}[Relation to bootstrap]
The idea to use the master loop equation to deduce results about Wilson loop expectations has also recently appeared in physics, under the name ``bootstrap". See the recent works \cite{kazakov2023bootstrap, Kazakov2024bootstrap} which obtain numerical results along these lines. In other areas of physics, the bootstrap has been quite successful in numerically determining observable expectation values, see e.g. the breakthrough work \cite{IsingBootstrap} on the 3D Ising model.
\end{remark}

\begin{remark}[Relation to the surface sums of \cite{CPS2023}]\label{remark:surface-exploration}
The work \cite{CPS2023} gave a formula for Wilson loop and string expectations in terms of weighted surface sums. Given a loop, the formula involves a sum over all surfaces spanning the loop. Although these surface sums do not appear directly in the present work, we found it conceptually quite helpful to think in terms of surfaces. In terms of surfaces, our proof strategy proceeds by exploring the surface spanned by the loop, one face at a time. This is exactly analogous to how one analyzes models of random planar maps via peeling \cite{curien2016planar}. In the context of planar maps, conditioning on the first step of the peeling process leads to Tutte's equation, which is absolutely central in the analysis of planar maps. In our context, conditioning on the first step of the exploration process of our surfaces precisely leads to the master loop equation. Thus, one should think of the master loop equation as the analog of Tutte's equation for the surface models described in \cite{CPS2023}.

In terms of the surfaces, it becomes clear why one should expect area law (at least in our regime of $\beta$): every time we see a plaquette in the exploration, this comes with a huge cost. Thus if we know that we must see at least $\area(\ell)$ many plaquettes (since our sum only involves surfaces which span the loop $\ell$), we should expect to have to pay a cost which is exponential in $\area(\ell)$. The way to make this rigorous is via the contraction estimate \eqref{eq:contraction}.

From the surface viewpoint, the problematic merger term is related to understanding the contribution of higher genus surfaces. On the one hand, each surface is weighted by their genus $g$ like $N^{-2g}$. Thus, when $N$ is large, we prefer to have lower genus surfaces. This gives another heuristic towards why it should be easier to prove area law at large $N$, because in general, sums over planar maps or low-genus surfaces are easier to understand than sums over very high genus surfaces. Nevertheless, when $N$ is large but finite, we still need to quantitatively understand the contributions of the very high genus surfaces, which in general is a difficult problem. See e.g. \cite{budzinski2021local, budzinski2022local} for recent results in this vein. 

Finally, for a surface interpretation of the extension argument from the truncated model to the original model, we can imagine a modified exploration process, where we allow ourselves only to explore those plaquettes which have a limited number of copies (recall the upper bound we imposed in the truncated model). As previously mentioned, it is extremely rare for a plaquette to have too many copies, and so if we ever come across such a plaquette, we already know there is a huge associated cost, and so there is no need to explore at this plaquette.
\end{remark}

We briefly summarize the remainder of the paper. In Section \ref{section:preliminaries}, we make definitions and introduce some preliminary results. In Section \ref{section:truncated-model}, we introduce the truncated model, and then show that area law holds for this model. In Section \ref{section:original-model}, we extend the argument to the original model, ultimately proving Theorem \ref{thm:area-law-main}. \\

\noindent \textbf{Acknowledgements:} We thank Jacopo Borga and Jasper Shogren-Knaak for helpful conversations. S.C. was partially supported by the NSF under Grant No. DMS-2303165. R.N. was supported by the NSF under Grant No. GRFP-2141064. S.S. was partially supported by the NSF under Grant No. DMS-2153742.

\section{Preliminaries}\label{section:preliminaries}

The present section serves as a preliminary section which will set notation and discuss preliminary results to be used in the rest of the paper.  First, we fix some notation. The parameter $N \geq 1$ denotes the size of the matrices in our matrix group $\unitary(N)$. In this paper, whenever we refer to a lattice $\Lambda$, we will mean a square lattice of the form $\Lambda=[-L,L]^d \sse \Z^d$. One should think of $L \gg 1$, and in this paper, all of our estimates will be uniform in $L$. Let $E_{\Lambda}^+$ (resp. $E_\Lambda$) be the collection of positively oriented (resp. oriented) edges in $\Lambda$. Similarly, let $\mc{P}_{\Lambda}^+$ (resp. $\mc{P}_{\Lambda}$) denote the set of positively oriented (resp. oriented) plaquettes in $\Lambda$. For oriented edges $e \in E_\Lambda$, we denote by $e^{-1}$ the oppositely oriented version of $e$. Similarly, for oriented plaquettes $p \in \mc{P}_\Lambda$, we denote by $p^{-1}$ the oppositely oriented version of $p$.

\begin{definition}[Lattice gauge configuration]
A lattice gauge configuration $U$ is a function $U : E_\Lambda^+ \ra \unitary(N)$. We always implicitly extend $U : E_\Lambda \ra \unitary(N)$ to oriented edges, by imposing that $U_{e^{-1}} := U_e^{-1}$ for all $e \in E_\Lambda^+$.
\end{definition}

\begin{definition}[Plaquette variable]
Given a lattice gauge configuration $U : E_\Lambda^+ \ra \unitary(N)$ and an oriented plaquette $p = e_1 e_2 e_3 e_4$, we define the plaquette variable (abusing notation)
\begin{equs}
U_p := U_{e_1} U_{e_2} U_{e_3} U_{e_4}.
\end{equs}
\end{definition}


Next, we define the lattice Yang--Mills measure.

\begin{definition}[Lattice Yang--Mills]\label{def:lattice-ym}
The \textit{lattice Yang-Mills measure} with Wilson action and gauge group $\mathrm{U}(N)$ is the measure on $\mathrm{U}(N)^{E_{\Lambda}^+}$ given by
\begin{equs}
    d\mu_{\Lambda, \beta, N}(U) = \frac{1}{Z_{\Lambda, \beta, N}} \prod_{p \in \mc{P}_{\Lambda}^+} \exp(2\beta \mathrm{Re}\mathrm{Tr}(U_p)) dU,
\end{equs}
where $dU = \prod_{e \in E_\Lambda^+} dU_e$ is an independent copy of Haar measure on $\mathrm{U}(N)$ for each edge $e \in E_{\Lambda}^+$, and 
\begin{equs}
Z_{\Lambda, \beta, N} :=\int_{\mathrm{U}(N)^{E_{\Lambda}^+}}\prod_{p \in \mc{P}_{\Lambda}^+} \exp(2\beta \mathrm{Re}\mathrm{Tr}(U_p)) dU.
\end{equs}
We denote by $\langle \cdot \rangle_{\Lambda, \beta, N}$ expectation with respect to the probability measure $\mu_{\Lambda, \beta, N}$.
\end{definition}

\begin{remark}[Trace versus normalized trace]\label{remark:trace-vs-normalized-trace}
In this paper, we choose to use the trace instead of the normalized trace in the definition of lattice Yang--Mills. On the other hand, the normalized trace is a natural choice from certain points of view. For instance, $\tr(U) = \frac{1}{N} \Tr(U) \in [-1, 1]$ for unitary matrices $U \in \unitary(N)$, for any $N \geq 1$. That is, the range of the normalized trace does not depend on $N$. Note also that $\beta \Tr(U) = \tilde \beta \tr(U)$, where $\tilde \beta = N \beta$ is the scaled analog of $\beta$ mentioned in the caption of Figure \ref{figure:os-comparison} (right). For $x \in [-1,1]$, one can approximate the function $x \to e^{2 \tilde \beta x}$ by a polynomial of $x$ that does not depend on $N$.
\end{remark}

Next, we define a collection of quantities associated to the lattice Yang--Mills measure. These quantities are associated to loops and strings, which we first define.

\begin{definition}[Loops and strings]\label{def:loops-and-strings}
We will represent loops $\ell$ in $\Lambda$ by the sequence of oriented edges $\ell = e_1 \cdots e_n$ that are traversed by $\ell$. We denote $|\ell| := n$ the perimeter (i.e., length) of the loop. For $x \in \Z^d$, we denote by $\varnothing_x$ the null loop at $x$, i.e. the loop starting and ending at $x$ with no edges.

We say that a loop $\ell$ contains a backtrack if $\ell$ is of the form $\ell = a ee^{-1} b$. We say that $\ell' = a b$ is the loop obtained by erasing the backtrack $ee^{-1}$. We denote by $[\ell]$ the loop obtained by erasing all backtracks in $\ell$ (note the order of erasure does not matter).

A \textit{string} $s=(\ell_1,\dots,\ell_n)$ is a sequence of loops $\ell_1,\dots,\ell_n$. We let $|s|:=\sum_{i=1}^{n}|\ell_i|$ denote the perimeter of the string $s$. Let $\mc{S}_\Lambda$ be the set of finite collections of nonempty loops $s = (s_1, \ldots, s_n)$ contained in $\Lambda$ with no backtracks. By convention, the null loops $\varnothing_x$ are elements of $\mc{S}_\Lambda$.
\end{definition}



\begin{definition}[Loop variables and observables]
Let $U : E_\Lambda^+ \ra \unitary(N)$ be a lattice gauge configuration. For  loops $\ell=e_1e_2\dots e_n$, let $U_{\ell}:=U_{e_1}U_{e_2}\dots U_{e_n}$ denote the corresponding loop variable. We define
\begin{equs}
W_\ell(U) := \tr(U_\ell),
\end{equs}
where $\tr = \frac{1}{N} \Tr$ is the normalized trace. For strings $s=(\ell_1,\dots,\ell_n)$, we define
\begin{equs}
W_s(U):= \prod_{i \in [n]} W_{\ell_i}(U).
\end{equs}
We refer to $W_\ell$ as a Wilson loop observable, and $W_s$ as a Wilson string observable.
\end{definition}

One can think of the Wilson string expectations as the multipoint functions of the lattice Yang--Mills measure, and one typically studies the measure via its multipoint functions. We remark that at this point, we have introduced all the notation appearing in Theorem \ref{thm:area-law-main}.

Next, due to the identity $UU^{-1} = \groupid$ for unitary matrices $U$, we have that $U_\ell = U_{[\ell]}$ and $W_\ell(U) = W_{[\ell]}(U)$ for all loops $\ell$, that is, the Wilson loop observables are invariant under backtrack erasure. Note also that $U_{\varnothing_x} = \groupid$ and $W_{\varnothing_x}(U) = 1$.

\begin{notation}
Often in this paper, we will have a function $J : \mc{P}_\Lambda \ra \N$, and we will raise a number $x$ to the power $J$, i.e. we will write
\begin{equs}
x^J := \prod_{p \in \mc{P}_\Lambda} x^{J(p)} = x^{\sum_{p \in \mc{P}_\Lambda} J(p)}.
\end{equs}
We will also write $B - J$ where $B$ is a number. This should be taken to mean the function $\mc{P}_\Lambda \ra \N$ which takes value $(B - J)(p) = B - J(p)$. The same considerations apply if $J : \mc{P}_\Lambda^+ \ra \N$ is a function of the positively oriented plaquettes.
\end{notation}


\subsection{Parameters}

Throughout this paper, we will assume that $N$ is large enough depending on $d$ such that
\begin{equs}\label{eq:N-large-assumption}
N \geq 10^{10} d^{10}.
\end{equs}
In Section \ref{section:original-model}, we will comment on the case where $N$ is smaller than $10^{10}d^{10}$. We will take parameters $\beta, B \geq 0$ such that
\begin{equs}
\beta &\leq 10^{-10d} d^{-1}, \label{eq:beta-choice} \\
~~ \frac{1}{2d} 10^{-3} N  \leq B &\leq \frac{1}{d} 10^{-3} N. \label{eq:B-choice}
\end{equs}
Moreover, we will assume that $B$ is odd, which will be convenient later. One should always think of relative sizes of $\beta, B$ as satisfying
\begin{equs}
d N \beta \ll B, ~~ d^5 \ll B, ~~ dB \ll N.
\end{equs}
In particular, at many intermediate places, there will be constants $C$ floating around, and due to our choices of $N, \beta, B$, we will always e.g. have that
\begin{equs}
C d N \beta \leq 10^{-3}B, ~~ Cd^5 \leq 10^{-3} B.
\end{equs}

\begin{remark}
In plain terms, one should think of $\beta$ as tiny, $B$ as huge, and $N$ as even larger than $B$ by a dimensional factor. Of course, the precise values such as $10^{10} d^{10}$ are somewhat arbitrary, and are chosen so that we always have a lot of room in the estimates.
\end{remark}

The assumptions \eqref{eq:N-large-assumption}-\eqref{eq:B-choice} are always implicitly assumed to hold (unless otherwise stated). Nonetheless, at various points in the statements of results we will still explicitly recall these parameter assumptions.

\begin{definition}[Truncated exponential]\label{def:truncated-exponential}
For $k_0 \geq 0$, define the truncated exponential 
\begin{equs}
\exp_{k_0}(x):=\sum_{k=0}^{k_0} \frac{x^k}{k!}.
\end{equs}
Define also the exponential tail function
\begin{equs}
\tail_{k_0} (x) := \exp(x) - \exp_{k_0}(x) = \sum_{k = k_0 + 1}^\infty \frac{x^k}{k!}.
\end{equs}
\end{definition}

\begin{remark}
As we will see later in Section \ref{section:truncated-model}, the truncated model will be defined by replacing the usual exponential $\exp$ in Definition \ref{def:lattice-ym} by the truncated exponential $\exp_B$. One outcome of the next lemma is that $\tail_B$ is tiny compared to $\exp_B$, which is related to the discussion in Section \ref{section:introduction} about why it is extremely rare for a plaquette to have too many copies.
\end{remark}

We will need the following properties of the truncated exponential.

\begin{lemma}[Properties of truncated exponential]\label{lemma:properties-of-truncated-exp}
Under the parameter assumptions \eqref{eq:N-large-assumption}-\eqref{eq:B-choice}, the following properties hold for all $|x| \leq 2N\beta$:
\begin{enumerate}
\item    \begin{equs}\label{eq:truncated<full}
    0 < \exp_B(x) \leq \exp(x).
\end{equs}
\item 
\begin{equs}\label{eq:Truncated vs Untruncated exp}
|\tail_B(x)| \leq \exp(-B/10) \exp_B(x).
\end{equs}
\item For all $k \leq B$,
\begin{equs}\label{eq:Truncation level comparison}
|\exp_k(x)| \leq e^{B-k}\exp_B(x).
\end{equs}
\end{enumerate}
\end{lemma}

\begin{proof}
We begin by proving one half of \eqref{eq:truncated<full}. Since $B$ is odd and $B \geq 4 (2N\beta)$, for all $|x| \leq 2N\beta$,
\begin{equs}
    \exp(x)-\exp_B(x)&=\frac{x^{B+1}}{(B+1)!}+\sum_{k=B+2}^{\infty} \frac{x^k}{k!}\geq \frac{x^{B+1}}{(B+1)!} \bigg(1- \sum_{k=1}^{\infty} \bigg(\frac{|x|}{B+2}\bigg)^k\bigg)\\
    &\geq \bigg(1-\sum_{k=1}^{\infty}(1/4)^k\bigg)\frac{x^{B+1}}{(B+1)!}>0,
\end{equs}
proving that $\exp(x)-\exp_B(x) \geq 0$. 

To prove the other half of \eqref{eq:truncated<full} as well as \eqref{eq:Truncated vs Untruncated exp}, we recall the classical Poisson concentration inequality (see e.g. \cite[Chapter 2]{BLM2013}):
\begin{equs}\label{eq:poisson-concentration}
\p(X \geq \lambda + t) \leq \exp\bigg(-\frac{t^2}{2(\lambda + t)}\bigg), ~~ t \geq 0, ~~ X \sim \mrm{Poisson}(\lambda).
\end{equs}
Applying this with $\lambda = 2N \beta$ and $t = B + 1 - \lambda  \gg \lambda$, we obtain
\begin{equs}
e^{-2N\beta} \sum_{k=B+1}^\infty \frac{(2N\beta)^k}{k!} \leq \exp(-B/4),
\end{equs}
which implies that for $|x| \leq 2N\beta$,
\begin{equs}\label{eq:tail-intermediate-upper-bd}
|\tail_B(x)| \leq e^{2N\beta} \exp(-B/4) \leq \exp(-B/6).
\end{equs}
It follows that
\begin{equs}
\exp_B(x) = \exp(x) - \tail_B(x) &\geq \exp(x) - |\tail_B(x)| \\
&\geq \exp(-2N\beta) - \exp(-B/6) \\
&\geq \exp(-2N\beta) \big(1 - \exp(-B/6 + 2N\beta)\big) \\
&\geq \exp(-2N\beta)(1 - \exp(-B/8)).
\end{equs}
In summary, we obtain the lower bound
\begin{equs}\label{eq:truncated-exp-lower-bound}
\exp_B(x) \geq \exp(-2N\beta)(1 - \exp(-B/8)) > 0,
\end{equs}
which proves the other half of \eqref{eq:truncated<full}. Next, to prove \eqref{eq:Truncated vs Untruncated exp}, we combine the above with \eqref{eq:tail-intermediate-upper-bd}, and claim that
\begin{equs}
\exp(-B/6) \leq \exp(-B/10) \exp(-2N\beta) (1 - \exp(-B/8)).
\end{equs}
This can be seen by rearranging and using our parameter assumptions, since $B \gg N \beta$ and $B \gg d \geq 1$. \\

\noindent Lastly, we prove \eqref{eq:Truncation level comparison} by splitting into two cases. \\

\noindent \emph{Case 1: $B/10 \leq k \leq B$.}
If $k = B$ there is nothing to prove. Thus, assume further that $k < B$. Again by Poisson concentration \eqref{eq:poisson-concentration}, we have that (since $k - 2N\beta \gg 2N\beta$)
\begin{equs}
e^{-2N\beta} |\exp_k(x) - \exp_B(x)| &\leq \exp(-(k-2N\beta)/4) + \exp(-(B - 2N\beta)/4) \\
&\leq  \exp(-B/50).
\end{equs}
It follows that
\begin{equs}
|\exp_k(x)| &\leq \exp(2N\beta - B/50) + e^{-2N\beta} \exp_B(x) \\
&\leq 2 \exp_B(x) \leq 2^{B-k} \exp_B(x)
\end{equs}
where we used the lower bound \eqref{eq:truncated-exp-lower-bound}, as well as the fact that $k < B$. \\

\noindent \emph{Case 2: $k < B/10$.}
In this case, we may bound
\[ e^{k-B} |\exp_k( x)| \leq e^{-9B/10} \sum_{j=0}^k \frac{|x|^k}{k!} \leq e^{-9B/10} \exp(2N \beta).  \]
To finish, we again use the lower bound \eqref{eq:truncated-exp-lower-bound} to reduce to the claim that
\begin{equs}
e^{-9B/10} \exp(2N\beta) \leq \exp(-2N\beta)(1 - \exp(-B/8)),
\end{equs}
which follows from our parameter choices.
\end{proof}

\subsection{General master loop equation}\label{section:mle}

In this section, we recall the master loop equation which will be central to our analysis. At a high level, the master loop equation relates a given string expectation to a weighted sum of string expectations, over strings which ``neighbor" the original string. Such neighboring strings are obtained by performing one of a finite set of string operations, which we next discuss.

Let $s = (s_1, \ldots, s_n) \in \mc{S}_\Lambda$, that is, $s$ is a string without backtracks. To define our string operations, will assume that $s$ is not a null loop (otherwise, there would be no operations). There are six total types of string operations, grouped into three categories: splittings, mergers, and deformations.

\begin{notation}
We reserve the letter $e$ for edges of the lattice $\Lambda$, and $\mbf{e}$ for edges of strings $s$.
\end{notation}

\begin{definition}[Splitting]

Suppose that $\ell_1$ contains two edges $\mbf{e}, \mbf{e}'$ which are the same lattice edge, so that $\ell_1 = a \mbf{e} b \mbf{e}' c$. Define loops $(\ell_{1, 1}', \ell_{1, 2}') := (a \mbf{e}' c, \mbf{e} b)$. The \textit{positive splitting} of $s$ at $(\mbf{e}, \mbf{e}')$ is defined to be the string $s' = ([\ell_{1, 1}'], [\ell_{1, 2}'], \ell_2, \ldots, \ell_n)$. 

Similarly, if $\ell_1$ contains edges $\mbf{e}, \mbf{e}'$ such that $\mbf{e}, (\mbf{e}')^{-1}$ are the same lattice edge, so that $\ell_1 = a \mbf{e} b \mbf{e}' c$ or $\ell_1 = a \mbf{e}' b \mbf{e} c$, then define loops $(\ell_{1, 1}', \ell_{1, 2}') := (ac, b)$. The \textit{negative splitting} of $s$ at $(\mbf{e}, \mbf{e}')$ is defined to be the string $s' := ([\ell_{1, 1}'], [\ell_{1, 2}'], \ell_2, \ldots, \ell_n)$.

The definitions directly extend to the case where the edges $\mbf{e}, \mbf{e}'$ are on some other loop $\ell_j$. Given an edge $\mbf{e}$ of $s$, we let $\mbb{S}_{+}(s, \mbf{e}) 
/ \mbb{S}_-(s, \mbf{e})$ to be the sets of all positive/negative splittings of $s$ at $(\mbf{e}, \mbf{e}')$ for some other edge $\mbf{e}'$ (which must necessarily be on the same loop as $\mbf{e}$). We will frequently use the notation $\mbb{S}_{\pm}(s, \mbf{e})$ to denote the set of positive or negative splittings of $s$ at $\mbf{e}$.

\end{definition}
     
\begin{definition}[Merger]
Let $\ell_i, \ell_j$ be loops of $s$ which both contain some lattice edge $e$. Thus, $\ell_i = a \mbf{e} b$ and $\ell_j = c \mbf{e}' d$, where $\mbf{e}, \mbf{e}'$ are the same lattice edge. Define $\ell_i \oplus_{\mbf{e}, \mbf{e}'} \ell_j := a \mbf{e}' d c \mbf{e} b$. The \textit{positive merger} of $s$ at $(\mbf{e}, \mbf{e}')$ is defined to be the string $s' := ([\ell_i \oplus_{\mbf{e}, \mbf{e}'} \ell_j], \ell_1, \ldots, \hat{\ell}_i, \ldots, \hat{\ell}_j, \ldots, \ell_n)$, where here $\hat{\ell}$ denotes omission.

Similarly, if $\ell_i = a \mbf{e} b$, $\ell_j = c \mbf{e}' d$ where $\mbf{e}, (\mbf{e}')^{-1}$ are the same lattice edge, then define $\ell_i \ominus_{\mbf{e}, \mbf{e}'} \ell_j := adcb$. The \textit{negative merger} of $s$ at $(\mbf{e}, \mbf{e}')$ is defined to be the string $s' := ([\ell_i \oplus_{\mbf{e}, \mbf{e}'} \ell_j], \ell_1, \ldots, \hat{\ell}_i, \ldots, \hat{\ell}_j, \ldots, \ell_n)$.

Given an edge $\mbf{e}$ of $s$, we let $\mbb{M}_{+}(s, \mbf{e})/\mbb{M}_-(s, \mbf{e})$ be the sets of all positive/negative mergers of $s$ at $(\mbf{e}, \mbf{e}')$ for some other edge $\mbf{e}'$. We will frequently use the notation $\mbb{M}_{\pm}(s, \mbf{e})$ to denote the set of positive or negative mergers of $s$ at $\mbf{e}$.

\end{definition}

\begin{definition}[Deformation]
Let $\mbf{e}$ be an edge of $\ell_i$ for some loop $\ell_i$ of $s$. Then $\ell_i = a \mbf{e} b$. For any oriented plaquette $p \in \mc{P}_\Lambda$ which contains (the lattice edge corresponding to) $\mbf{e}$, we may consider the loop given by $\ell_p = \mbf{e}' c$ tracing the boundary of $p$. Here, $\mbf{e}'$ is the same lattice edge as $\mbf{e}$. We define $\ell_i \oplus_\mbf{e} p := a \mbf{e}' c \mbf{e} b$, and define the \textit{positive deformation} of $s$ at $\mbf{e}$ to be the string $s' := (\ell_1, \ldots, [\ell_i \oplus_{\mbf{e}} p], \ldots, \ell_n)$.

Similarly, for any oriented plaquette $p \in \mc{P}_\Lambda$ which contains (the lattice edge corresponding to) $\mbf{e}^{-1}$, we may consider the loop given by $\ell_p = \mbf{e}' c$ tracing the boundary of $p$. Here, $\mbf{e}'$ is the same lattice edge as $\mbf{e}^{-1}$. We define $\ell_i \ominus_\mbf{e} p := a c b$, and define the \textit{negative deformation} of $s$ at $\mbf{e}$ to be the string $s' := (\ell_1, \ldots, [\ell_i \oplus_{\mbf{e}'} p], \ldots, \ell_n)$

Given an edge $\mbf{e}$ of $s$, we let $\mbb{D}_{+}(s, \mbf{e})/ \mbb{D}_-(s, \mbf{e})$ to be the set of all positive/negative deformations of $s$ at $\mbf{e}$. We will frequently use the notation $\mbb{D}_{\pm}(s, \mbf{e})$ to denote the set of positive or negative deformations of $s$ at $\mbf{e}$.

In principle, the definition of $\mbb{D}_{\pm}$ depends on $\Lambda$, since we restrict to only considering plaquettes $p \in \mc{P}_\Lambda$. We will omit this dependence in the following, to lessen notation. 
\end{definition}

\begin{remark}[Reduced loops may vanish]
In some cases (but not all), the loop(s) obtained after applying a string operation may actually reduce to zero. For instance, if the first loop of $s = (\ell_1, \ell_2)$ is a single plaquette loop corresponding to some plaquette $p$, and we negatively deform using the oppositely oriented plaquette $p^{-1}$, then reduced loop $[\ell_1 \ominus_{\mbf{e}} p^{-1}]$ will actually be a null loop. In such situations, we erase any null loops arising in $s'$, so that $s' = (\ell_2)$ instead of $(\varnothing_x, \ell_2)$. With this convention, all string operations result in elements of $\mc{S}_\Lambda$, which recall was defined (in Definition \ref{def:loops-and-strings}) to be the set of strings $s$ which are either null loops, or such that all loops of $s$ are nonempty and without backtracks.
\end{remark}


    


Next, we prove a general master loop equation which will be used in Sections \ref{section:truncated-model} and \ref{section:original-model}.

\begin{prop}[General master loop equation at a single location]\label{prop:general-mle}
Let $\varep > 0$ and let $\mbb{D}(0, 2N\beta + \varep)$ denote the open ball in $\C$ of radius $2N\beta + \varep$ centered at $0$. For each $p \in \mc{P}_\Lambda^+$, fix a holomorphic function $\rho_p: \mathbb{D}(0,2N\beta+\varep) \to \C$. Then for any non-null string $s$ and any edge $\mbf{e}$ of $s$, we have that
\begin{equs}
\int dU W_s(U)\prod_{p \in \mc{P}_{\Lambda}^{+}} &\rho_p(2\beta \mathrm{Re}\mathrm{Tr}(U_p)) \\
&= \mp \sum_{s' \in \mbb{S}_{\pm}(s, \mbf{e})} \int dU W_{s'}(U)\prod_{p \in \mc{P}_{\Lambda}^{+}} \rho_p(2\beta \mathrm{Re}\mathrm{Tr}(U_p)) \\
&\mp \frac{1}{N^2} \sum_{s' \in \mbb{M}_{\pm}(s, \mbf{e})} \int dU W_{s'}(U)\prod_{p \in \mc{P}_{\Lambda}^{+}} \rho_p(2\beta \mathrm{Re}\mathrm{Tr}(U_p)) \\
&- \frac{\beta}{N} \sum_{\substack{s' \in \mbb{D}_{+}(s, \mbf{e})\\
s' = s \oplus_{\mbf{e}}p_0 }} \int dU W_{s'}(U)\rho_{p_0}'(2\beta \mathrm{Re}\mathrm{Tr}(U_{\tilde{p}}))\prod_{p \in \mc{P}_{\Lambda}^{+} \backslash \{p_0\}} \rho_p(2\beta \mathrm{Re}\mathrm{Tr}(U_p)) \\
&+ \frac{\beta}{N} \sum_{\substack{s' \in \mbb{D}_{-}(s, \mbf{e})\\
s' = s \ominus_{\mbf{e}}p_0 }} \int dU W_{s'}(U)\rho_{p_0}'(2\beta \mathrm{Re}\mathrm{Tr}(U_{\tilde{p}}))\prod_{p \in \mc{P}_{\Lambda}^{+} \backslash \{p_0\}} \rho_p(2\beta \mathrm{Re}\mathrm{Tr}(U_p)).
\end{equs}
To be clear, the $\mp$ notation is paired with the $\pm$ notation, so that positive splittings and mergers come with a minus sign, while negative splittings and mergers come with a positive sign. We refrain from using the notation for deformations, to avoid introducing an ``oplusminus" command in \LaTeX. Also, in the last two terms, for $s' \in \mbb{D}_{\pm}(s, \mbf{e})$, we write $p_0$ for the plaquette which is used to deform $s$.
\end{prop}

\begin{proof}
Since $\rho$ is holomorphic on $\mathbb{D}(0,2N\beta+\epsilon)$, it has a Taylor series centered at $z=0$ which converges absolutely on $\mathbb{D}(0,2N\beta)$. We thus have that
\begin{equation}\label{eq:integral-taylor-expansion}
\begin{aligned}
\int &dU W_s(U)\prod_{p \in \mc{P}_{\Lambda}^{+}} \rho_p(2\beta \mathrm{Re}\mathrm{Tr}(U_p))\\
&=\sum_{K: \mc{P}_{\Lambda}^+ \to \N} \frac{\beta^{K}}{K!}\int dU W_s(U)\prod_{p \in \mc{P}_{\Lambda}^+} \rho_p^{(K(p))}(0) (2\mathrm{Re}\mathrm{Tr}(U_p))^{K(p)}\\
&=\sum_{J: \mc{P}_{\Lambda} \to \N} \frac{\beta^{J}}{J!} \prod_{p \in \mc{P}_\Lambda^+} \rho_p^{(K_J(p))}(0) \int dU W_s(U)\prod_{p \in \mc{P}_{\Lambda}} \mathrm{Tr}(U_p)^{J(p)},
\end{aligned}
\end{equation}
where $K_J(p) := J(p) + J(p^{-1})$ for $p \in \mc{P}_\Lambda^+$. Here, in the last identity we used that $2\mrm{Re}\Tr(U_p) = \Tr(U_p) + \Tr(U_{p^{-1}})$ as well the binomial expansion identity:
\begin{equs}
\frac{1}{k!} \sum_{j=0}^k (x + y)^k = \sum_{j=0}^k \frac{1}{j!(k-j)!} x^j y^{k-j}.
\end{equs}
The argument now follows similarly to  \cite[Theorem 5.7]{CPS2023}. For brevity, let
\[
I(s,J) := \int dU W_s(U) \prod_{p \in \mathcal{P}} \operatorname{Tr}(U_p)^{J(p)} .
\]
Fix $J : \mathcal{P} \to \mathbb{N}$. It may help to keep in mind that $J(p)$ counts the number of copies of $p$ that are present. Let us now set some notation. Let $e$ be the lattice edge given by 
$\mbf{e}$. We write $p \succ e$ to mean that $p$ contains $e$, and $p \succ {e^{-1}}$ to means that $p$ contains $e^{-1}$. If $p \succ e$ or $p \succ {e^{-1}}$, let $s \oplus_{\mbf{e}} p$ and $s \ominus_{\mbf{e}} p$ be the positive and negative deformations of $s$ by $p$ at $\mbf{e}$. For $p \in \mathcal{P}_\Lambda$, let $\delta_p : \mathcal{P}_\Lambda \to \mathbb{N}$ be the delta function at $p$. Following the proof of \cite{CPS2023}[Theorem 5.7], by applying the word recursion \cite[Proposition 5.3]{CPS2023}, we obtain that
\[
\begin{aligned}
I(s,J) = & - \sum_{s' \in \mathbb{S}_+(s, \mbf{e})} I(s',J) + \sum_{s' \in \mathbb{S}_-(s, \mbf{e})} I(s',J) \\
& - \frac{1}{N^2} \sum_{s' \in \mathbb{M}_+ (s, \mbf{e})} I(s',J) + \frac{1}{N^2} \sum_{s' \in \mathbb{M}_-(s, \mbf{e})} I(s',J) \\
& - \frac{1}{N} \sum_{\substack{p \in \mathcal{P} \\ p \succ e}} J(p) I(s \oplus_{\mbf{e}} p, J - \delta_p) 
+ \frac{1}{N} \sum_{\substack{p \in \mathcal{P} \\ p \succ {e^{-1}}}} J(p) I(s \ominus_{\mbf{e}} p, J - \delta_p).
\end{aligned}
\]
(Here, the factor of $J(p)$ arising in the last two terms arises because there are $J(p)$ copies of the plaquette $p$ which can possibly be used to deform $s$.) Combining this with the identity \eqref{eq:integral-taylor-expansion}, we obtain
\[
\begin{aligned}
\int dU W_s(U)\prod_{p \in \mc{P}_{\Lambda}^{+}} &\rho_p(2\beta \mathrm{Re}\mathrm{Tr}(U_p)) \\
&=  \mp\sum_{s' \in \mathbb{S}_\pm{(s, \mbf{e})}} \int dU W_{s'}(U)\prod_{p \in \mc{P}_{\Lambda}^{+}} \rho_p(2\beta \mathrm{Re}\mathrm{Tr}(U_p))\\
&\mp\frac{1}{N^2} \sum_{s' \in \mathbb{M}(s, \mbf{e})} \int dU W_{s'}(U)\prod_{p \in \mc{P}_{\Lambda}^{+}} \rho_p(2\beta \mathrm{Re}\mathrm{Tr}(U_p)) \\
& + D_1 + D_2,
\end{aligned}
\]

where
\[
\begin{aligned}
D_1 & := - \frac{1}{N} \sum_{\substack{p_0 \in \mathcal{P} \\ p_0 \succ e }} \sum_{\substack{J : \mc{P}_\Lambda \ra \N \\ J(p_0) \geq 1}} 
J(p_0) \frac{\beta^J}{J!}  \bigg(\prod_{p \in \mc{P}_\Lambda^+} \rho_p^{(K_J(p))}(0)\bigg) I(s \oplus_{\mbf{e}} p_0, J - \delta_{p_0}), \\
D_2 & :=  \frac{1}{N} \sum_{\substack{p_0 \in \mathcal{P} \\ p_0 \succ {e^{-1}}}} \sum_{\substack{J : \mc{P}_\Lambda \ra \N \\ J(p_0) \geq 1}} J(p_0) \frac{\beta^J}{J!}   
\bigg(\prod_{p \in \mc{P}_{\Lambda}} \rho_p^{(K_J(p))}(0)\bigg) I(s \ominus_{\mbf{e}} p_0, J - \delta_{p_0}).
\end{aligned}
\]
Next, we change variables $J \mapsto J - \delta_{p_0}$ to obtain
\begin{equs}
D_1 = -\frac{\beta}{N} \sum_{\substack{p_0 \in \mathcal{P} \\ p_0 \succ e }} \sum_{\substack{J : \mc{P}_\Lambda \ra \N}} 
\frac{\beta^J}{J!} \rho_{p_0}^{(K_J(p_0) + 1)}(0) \bigg(\prod_{p \in \mc{P}_\Lambda^+ \backslash \{p_0\}} \rho_p^{(K_J(p))}(0)\bigg) I(s \oplus_{\mbf{e}} p_0, J).
\end{equs}
Finally, using the identity \eqref{eq:integral-taylor-expansion} in reverse, with $\rho_{p_0}$ replaced by $\rho_{p_0}'$, we obtain
\begin{equs}
D_1 = - \frac{\beta}{N} \sum_{\substack{s' \in \mbb{D}_{+}(s, \mbf{e})\\
s' = s \oplus_{\mbf{e}} p_0 }} \int dU W_{s'}(U)\rho_{p_0}'(2\beta \mathrm{Re}\mathrm{Tr}(U_{p_0}))\prod_{p \in \mc{P}_{\Lambda}^{+} \backslash \{p_0\}} \rho_p(2\beta \mathrm{Re}\mathrm{Tr}(U_p)).
\end{equs}
The term $D_2$ is handled similarly; we omit the details.
\end{proof}

\section{Truncated model}\label{section:truncated-model}

In this section, we define the following truncated version of the lattice Yang--Mills measure, which was previously alluded to in Section \ref{section:introduction}. We then proceed to prove area law for this truncated model. In the following, recall our parameter choices \eqref{eq:N-large-assumption}-\eqref{eq:B-choice}, as well as the definition of the truncated exponential (Definition \ref{def:truncated-exponential}).

\begin{definition}[Truncated lattice Yang--Mills]
Define the truncated lattice Yang--Mills measure via the density:
\begin{equs}
d\mu_{\Lambda, \beta, N, B} (U) := \frac{1}{Z_{\Lambda, \beta, N, B}} \prod_{p \in \mc{P}_\Lambda^+} \exp_B(2\beta \mrm{Re}\Tr(U_p)) dU.
\end{equs}
Let $\langle \cdot \rangle_{\Lambda, \beta, N, B}$ denote expectation with respect to the measure $\mu_{\Lambda, \beta, N, B}$.
\end{definition}

Since we are assuming throughout that our parameter choices \eqref{eq:N-large-assumption}-\eqref{eq:B-choice} are satisfied, we have by Lemma \ref{lemma:properties-of-truncated-exp} that the truncated exponential $\exp_B$ is pointwise positive, and so that $\mu_{\Lambda, \beta, N, B}$ is indeed a probability measure (in particular, it is positive). We now state the main theorem of the present section, which is that area law holds for the truncated model. 

\begin{theorem}[Area law for the truncated model]\label{thm:truncated-model-area-law}
With the choice of parameters as in \eqref{eq:N-large-assumption}-\eqref{eq:B-choice}, we have that for any rectangular loop $\ell$ contained in a lattice $\Lambda$,
\[ |\langle W_\ell \rangle_{\Lambda, \beta, N, B}| \leq 2  N^{|\ell|/4-1} (10^3 d \beta)^{\area(\ell)}. \]
In particular, we have area law decay for Wilson loop expectations in the truncated model.
\end{theorem}

\begin{remark}
Strictly speaking, Theorem \ref{thm:truncated-model-area-law} is not needed in the proof of Theorem \ref{thm:area-law-main}. However, its proof is quite short, and the proof of Theorem \ref{thm:area-law-main} in Section \ref{section:original-model} builds quite heavily on the ideas introduced in this section. Thus we chose to state and prove Theorem \ref{thm:truncated-model-area-law} in this section, as a warmup towards the full result Theorem \ref{thm:area-law-main}.
\end{remark}

To recall the discussion after the statement of Theorem \ref{thm:area-law-main}, the main idea is use the fact that Wilson string expectations satisfy a linear inhomogeneous equation, and then to prove a contraction estimate of the general form \eqref{eq:contraction}, which will result in an a priori bound on solutions to the equation. We begin to develop the necessary framework for the contraction mapping argument.

First, we define a notion of area for general strings $s$. In terms of the relation to the surface sums from \cite{CPS2023}, there is a very natural interpretation for the area we will define: it will precisely be the area of the minimal area surface which spans the string $s$. However, here we will directly define the area without introducing the surfaces appearing in \cite{CPS2023}.

\begin{definition}[Plaquette count]
We will refer to functions $K : \mc{P}_\Lambda^+ \ra \N$ or $J : \mc{P}_\Lambda \ra \N$ as plaquette counts. Our convention will be that the letter $K$ will denote an upper bound on the number of copies of each plaquette in our system, in that for each $p \in \mc{P}_\Lambda^+$, $K(p)$ is an upper bound on the total number of copies of $p$ and $p^{-1}$. Whereas when we use the letter $J$, we will think of $J(p)$ as specifying the exact number of copies of each oriented plaquette $p \in \mc{P}_\Lambda$. 

Given plaquette counts $K : \mc{P}_\Lambda^+ \ra \N$ or $J : \mc{P}_\Lambda \ra \N$, we define $\mrm{supp}(K) := \{p \in\mc{P}_\Lambda^+ : K(p) \neq 0\}$, and similarly for $\mrm{supp}(J)$. We refer to these sets as the supports of $K, J$.

Additionally, we will often abuse notation and write $B$ for the plaquette count $\mc{P}_\Lambda^+ \ra \N$ which assigns every plaquette $p \in \mc{P}_\Lambda^+$ the value $B$. It should always be clear from context whether $B$ represents a number or a plaquette count.
\end{definition}

\begin{definition}\label{def:balanced}
Let $s$ be a string, and let $J : \mc{P}_\Lambda \ra \N$ be a plaquette count. For oriented edges $e \in E_\Lambda$, define the edge occurrence number
\begin{equs}
n_e(s, J):= \text{number of occurrences of $e$ in $(s, J)$.}
\end{equs}
By this, we mean that each plaquette $p$ appears with $J(p)$ number of copies. Thus if a given edge $e$ is contained in $p$, then there will be $J(p)$ occurrences of $e$ coming from $p$. We say that $(s, J)$ is balanced if $n_e(s, J) = n_{e^{-1}}(s, J)$ for all $e \in E_\Lambda$.

Similarly, given a plaquette count $K : \mc{P}_\Lambda^+ \ra \N$ and a positively oriented edge $e \in E_\Lambda^+$, we define the ``unoriented" edge occurrence number
\begin{equs}
n_e(s, K) := \text{number of occurrences of $e$ or $e^{-1}$ in $(s, K)$.}
\end{equs}
In this definition, the convention is each plaquette $p \in \mc{P}_\Lambda^+$ has $K(p)$ copies.
\end{definition}

\begin{definition}[Area]\label{def:area}
Let $s$ be a string. We define $\area(s)$ to be the minimum area over all plaquette counts $J : \mc{P}_\Lambda \ra \N$ such that $(s, J)$ is balanced. I.e.,
\[ \area(s) := \min\Big\{\sum_{p \in \mc{P}_\Lambda} J(p) : (s, J) \text{ is balanced} \Big\}.\]
\end{definition}

\begin{remark}
Note that for rectangular loops $\ell$ with side lengths $R, L$, we have that $\area(\ell) = RL$, which coincides with the usual notion of area of a rectangular loop. This is because the minimal plaquette count is given by the indicator function on the set of plaquettes in the rectangle spanned by $\ell$.
\end{remark}

Having defined the area of a general string, next we define the quantities entering into the master loop equation. For $K : \mc{P}_\Lambda^+ \ra \N$, let
\begin{equs}\label{eq:phi-s-K}
\phi(s, K) := \int dU  W_s(U) \prod_{p \in \mc{P}_\Lambda^+} \exp_{K(p)}(2\beta \mrm{Re}\mrm{Tr}(U_p)) . \end{equs}
For null loops $\varnothing_x$, the value of $\phi(\varnothing_x, K)$ does not depend on the basepoint $x$, and is explicitly given by
\begin{equs}
\phi(\varnothing_x, K) = \int dU  \prod_{p \in \mc{P}_\Lambda^+} \exp_{K(p)}(2\beta \mrm{Re}\mrm{Tr}(U_p)). 
\end{equs}
Thus in a slight abuse of notation, we will write $\phi(\varnothing, K)$ to denote the above value. Also, note that
\begin{equs}\label{eq:partition-fn-identity}
\phi(\varnothing, B) = Z_{\Lambda, \beta, N, B}, \quad \langle W_s \rangle_{\Lambda, \beta, N, B} = \frac{\phi(s, B)}{Z_{\Lambda, \beta, N, B}}.
\end{equs}
Next, recall the string operations defined in Section \ref{section:mle}. We will need a variant of the deformation operation which takes the plaquette count $K$ into account.

\begin{definition}[Deformation]
Let $s \in \mc{S}_\Lambda$ be a string, and suppose that $s$ is non-null. Let $\mbf{e}$ be an edge of $s$. Let $K : \mc{P}_\Lambda \ra \N$ be a plaquette count. Let $\mbb{D}_{\pm}(s, \mbf{e}, K)$ be the set of pairs $(s', K')$ such that $s'$ is a positive/negative deformation of $s$ at $\mbf{e}$, and such that $K'$ satisfies the following constraint. First, let $p_0 \in \mc{P}_\Lambda^+$ be the positively oriented plaquette such that $s'$ is a deformation of $s$ using either $p_0$ or $p_0^{-1}$. Then we assume that $K(p_0) \geq 1$, and we impose that $K'$ is given by
\begin{equs}
K'(p) = \begin{cases} K(p) & p \neq p_0, \\
K(p) - 1 & p = p_0 .\end{cases}
\end{equs}
Note by definition that for all $(s', K') \in \mbb{D}_{\pm}(s, \mbf{e}, K)$, 
\begin{equs}\label{eq:defomration-area-change}
\sum_{p \in \mc{P}_\Lambda^+} K'(p) = -1 + \sum_{p \in \mc{P}_\Lambda^+} K(p).
\end{equs}
In the case where $K(p_0) = 0$ for all plaquettes $p_0$ such that either $p_0$ or $p_0^{-1}$ contain the edge $\mbf{e}$, then we define $\mbb{D}_{\pm}(s, \mbf{e}, K) := \varnothing$.
\end{definition}

\begin{remark}
Intuitively, one should think of $K$ as keeping track of the maximum number of plaquettes in the system, so that if we deform using a given plaquette, then the value of $K$ at that plaquette should be decremented.
\end{remark}

We next show the following master loop equation for $\phi(s, K)$. 

\begin{prop}[Master loop equation]\label{prop:mle}
Let $s \in \mc{S}_\Lambda$ be non-null, $\mbf{e}$ be an edge of $s$, and $K : \mc{P}_\Lambda^+ \ra \N$ be a plaquette count. We have that
\begin{equs}
\phi(s, K) = \mp \sum_{s' \in \mbb{S}_{\pm}(s, \mbf{e})} \phi(s', K) \mp \frac{1}{N^2} \sum_{s' \in \mbb{M}_{\pm}(s, \mbf{e})} \phi(s', K) \mp \frac{\beta}{N} \sum_{(s', K') \in \mbb{D}_{\pm}(s, \mbf{e}, K)} \phi(s', K').
\end{equs}
\end{prop}
\begin{proof}
This follows directly from Proposition \ref{prop:general-mle}, upon taking $\rho_p = \exp_{K(p)}$, and noting that $\exp_k' = \exp_{k-1}$, as well as that if $K(p) = 0$, then $\rho_p'(0) = 0$, so that in the deformation terms, one only has to restrict to those plaquettes for which $K(p) \geq 1$.
\end{proof}

Having obtained the master loop equation, we now turn towards defining the norm to be used in the contraction mapping argument. One key component of the norm involves the following statistic for strings, first introduced in \cite{Chatterjee2019a}, which we refer to as the ``splitting complexity". One should think of this statistic as keeping track of the number of possible splittings of any given string. In particular, we will show that this statistic decreases upon splitting. This will be essential in the ensuing contraction argument.

\begin{definition}[Splitting complexity]\label{def:splitting-complexity}
For non-null $s \in \mc{S}_\Lambda$, $s = (s_1, \ldots, s_n)$, let $\iota(s) := \frac{|s|}{4} - n$. Here, recall $|s|$ is the perimeter of $s$. For null loops $\varnothing_x$, we set $\iota(\varnothing_x) := 0$.
\end{definition}

In the following sequence of lemmas, we show how the splitting complexity changes under the various loop operations.

\begin{lemma}[Splitting complexity decreases after splitting]\label{lemma:splitting-complexity-after-splitting}
For $s \in \mc{S}_\Lambda$ non-null, an edge $\mbf{e}$ in $s$, and $s' \in \mbb{S}_{\pm}(s, \mbf{e})$, we have that $\iota(s') \leq \iota(s) - 1$. 
\end{lemma}
\begin{proof}
We first handle the positive splitting case, i.e. suppose that $s' \in \mbb{S}_+(s, \mbf{e})$. Without loss of generality, suppose that $\mbf{e}$ is the first edge of $s_1$, and so we may write $s_1 = \mbf{e} s_{1, 1} \mbf{e}' s_{1, 2}$. By the assumption that $s_1$ has no backtracks, neither do $s_{1, 1}, s_{1, 2}$. Thus neither $s_{1, 1}, s_{1, 2}$ are equivalent to the null loop, and thus we have that $s' = (\mbf{e}' s_{1, 1}, \mbf{e} s_{1, 2}, s_2, \ldots s_n)$. This shows that $\iota(s') = |s|/4 - (n+1) = \iota(s) - 1$.

Next, consider a negative splitting, i.e. $s' \in \mbb{S}_-(s, \mbf{e})$. Again, without loss of generality, suppose that $\mbf{e}$ is the first edge of $s_1$. Then $s_1 = \mbf{e} s_{1, 1} \mbf{e}' s_{1, 2}$
, where $\mbf{e}'$ is the same lattice edge as $\mbf{e}^{-1}$. Moreover, neither $s_{1, 1}, s_{1, 2}$ are empty (otherwise $s_1$ would have a backtrack) nor have backtracks. Thus $s' = (s_{1, 1}, s_{1, 2}, s_2, \ldots, s_n)$. Thus $\iota(s') = (|s|-2)/4 - (n+1) = \iota(s) - 3/2 \leq \iota(s) - 1$.
\end{proof}

\begin{lemma}[Splitting complexity increases by at most 1 after merger]\label{lemma:splitting-complexity-merger}
For $s \in \mc{S}_\Lambda$ non-null, an edge $\mbf{e}$ in $s$, and $s' \in \mbb{M}_{\pm}(s, \mbf{e})$, we have that $\iota(s') \leq \iota(s) + 1$.
\end{lemma}
\begin{proof}
Without loss of generality, suppose that $\mbf{e}$ is the first edge of the first loop $s_1$. First, suppose that $s'$ is a positive merger, i.e. $s' \in \mbb{M}_+(s, \mbf{e})$. Then $s_1 = \mbf{e} s_1'$, and for some $2 \leq j \leq n$,  $s_j = \mbf{e}' s_j'$, and we have that $s' = (\mbf{e}' s_j' \mbf{e} s_1', s_2, \ldots, \hat{s}_j, \ldots, s_n)$. Here, the $\hat{s}_j$ denotes omission. (Note that since $\mbf{e} s_1', \mbf{e}' s_j'$ do not have backtracks, neither does $\mbf{e}' s_1' \mbf{e} s_j'$.) Thus $\iota(s') = |s|/4 - (n-1) = \iota(s) + 1$.

Next, suppose $s'$ is a negative merger, i.e. $s' \in \mbb{M}_-(s, \mbf{e})$. Then $s_1 = \mbf{e} s_1'$, and for some $2 \leq j \leq n$,  $s_j = \mbf{e}' s_j'$, and we have that either $s' = ([s_1' s_j'], s_2, \ldots, \hat{s}_j, \ldots, s_n)$, or $s' = (s_2, \ldots, \hat{s}_j, \ldots, s_n)$, where the latter case occurs if $s_1' s_j'$ is equivalent to the null loop (recall that the notation $[s_1' s_j']$ denotes the loop $s_1' s_j'$ with any backtracks erased). In the former case, we have that $\iota(s') \leq (|s| - 2)/4 - (n-1) \leq \iota(s) + 1$. In the latter case, note that $s_1, s_j$ have length at least four (as they are closed loops in the lattice with no backtracks), and thus $|s'|\leq |s| - 8$, so that $\iota(s') \leq (|s|-8)/4 - (n-2) = |s| - n  = \iota(s)$.
\end{proof}

\begin{lemma}[Splitting complexity increases by at most 1 after deformation]\label{lemma:splitting-complexity-deformation}
Let $s \in \mc{S}_\Lambda$ be non-null, $\mbf{e}$ be an edge of $s$, $K : \mc{P}_\Lambda^+ \ra \N$ be a plaquette count, and $(s', K') \in \mbb{D}_{\pm}(s, \mbf{e}, K)$. 
Then $\iota(s') \leq \iota(s) + 1$.
\end{lemma}
\begin{proof}
Suppose without loss of generality that $\mbf{e}$ is the first edge of the first loop $s_1$. First suppose that $s'$ is a positive deformation, i.e. $(s', K') \in \mbb{D}_+(s, \mbf{e}, K)$. Then $s' = (s_1 \oplus_{\mbf{e}} p, s_2, \ldots, s_n)$ for some plaquette $p$ (there are no backtracks in $s_1 \oplus_{\mbf{e}} p$). So $4$ edges were added, thus $\iota(s') \leq (|s|+4)/4 - n = \iota(s) + 1$. Next, suppose that $s'$ is a negative deformation, i.e. $(s', K') \in \mbb{D}_-(s, \mbf{e}, K)$. Then either $s' = ([s_1 \ominus_{\mbf{e}} p], s_2, \ldots, s_n)$ (which happens in the case where $s_1 \ominus_{\mbf{e}} p$ is non-null, or $s' = (s_2, \ldots, s_n)$, which happens in the case where $s_1 \ominus_\mbf{e} p$ is null. In the former case, we have at most $3$ new edges, so $\iota(s') \leq (|s| + 3) / 4 - n \leq \iota(s) + 1$. In the latter case, since $|s_1| \geq 4$, we have that $\iota(s') \leq (|s| - 4)/4 - (n-1) = \iota(s)$.
\end{proof}

Next, we show how the area changes after applying a loop operation.

\begin{lemma}\label{lemma:area-after-loop-operation}
Let $s \in \mc{S}_\Lambda$ be non-null, and $\mbf{e}$ be an edge of $s$. For all $s' \in \mbb{S}_{\pm}(s, \mbf{e}) \cup \mbb{M}_{\pm}(s, \mbf{e})$, we have that $\area(s') = \area(s)$. For $s' \in \mbb{D}_{\pm}(s, \mbf{e})$, we have that $\area(s) \leq \area(s') + 1$.
\end{lemma}
\begin{proof}
From the definition of $\area$ as well as the definitions of splittings and mergers, we have that for $s' \in \mbb{S}_{\pm}(s, \mbf{e}) \cup \mbb{M}_{\pm}(s, \mbf{e})$ and any plaquette count $J$, the pair $(s, J)$ is balanced if and only if $(s', J)$ is balanced. In the case of deformations, suppose $s'$ is a deformation of $s$ using a plaquette $p$. Then for any plaquette count $J'$ such that $(s', J')$ is balanced, we have that $(s, J' + \delta_{p^{-1}})$ is balanced. Taking $J'$ to be minimal area, we obtain the inequality $\area(s) \leq \area(s') + 1$.
\end{proof}

Next, we define the space on which we will perform the contraction argument, as well as the norm that we will use. In the following, recall the unoriented edge occurrence number from Definition \ref{def:balanced}.

\begin{definition}[Truncated space]\label{def:truncated-space}
Define
\[ \Omega_B := \Big\{(s, K) : s \in \mc{S}_\Lambda , K \leq B, \max_{e \in E_\Lambda^+} n_e(s, K) \leq 2dB \Big\}. \]
Note that $\Omega_B$ is finite.
\end{definition}

\begin{remark}\label{remark:configuration-space}
In the definition of $\Omega_B$, we impose that $n_e(s, K) \leq 2dB$ for all $e \in E_\Lambda^+$. In words, this limits the number of copies of any given edge in our system. As we will see, this will impose an upper bound on the number of mergers (Lemma \ref{lemma:bound-on-loop-operations}), which is crucial for the contraction mapping argument.

We also remark that if $\ell$ is a rectangular loop, then $(\ell, B) \in \Omega_B$. To see this, note that each edge is contained in $2(d-1)$ unoriented plaquettes, and thus if each plaquette has at most $B$ copies, the number of copies of any given edge in $(\ell, B)$ is at most $1 + 2(d-1)B \leq 2dB$.
\end{remark}

\begin{definition}[Norm]\label{def:truncated-norm}
For $\lambda, \rho \leq 1$, and $\gamma \geq 1$, define the norm $\|\cdot\|_{\lambda, \gamma, \rho}$ on $\R^{\Omega_B}$ by
\[ \|f\|_{\lambda, \gamma, \rho} := \sup_{(s, K) \in \Omega_B} \lambda^{\iota(s)} \gamma^{\area(s)} \rho^{B-K} |f(s, K)|.\]
Here, $\rho^{B-K} = \prod_{p \in \mc{P}_\Lambda^+} \rho^{B - K(p)}$. Note that since $\Omega_B$ is finite, $\|f\|_{\lambda, \gamma, \rho} < \infty$ for any $f \in \R^{\Omega_B}$.
\end{definition}

\begin{remark}[Conceptual remarks]
On a conceptual level, one should think of $\lambda \ll 1$, so that the weight $\lambda^{\iota(s)}$ is tiny if $s$ has many possible splittings. The necessity of this weight is because in the master loop equation (Proposition \ref{prop:mle}), there is no small prefactor in front of the splitting term, while both the merger and deformation terms come with small prefactors. Chatterjee \cite{Chatterjee2019a} observed that by adding in this weight, a prefactor $\lambda$ will be generated in front of each splitting term, essentially because $\iota(s)$ decreases by at least 1 after any splitting. This small prefactor $\lambda$ then allows us to close the contraction argument.

On the other hand, in the merger and deformation terms, there will now be a large prefactor $\lambda^{-1}$, because  $\iota(s)$ may increase by 1 after a merger. Fortunately, the prefactors $N^{-2}$ and $\beta N^{-1}$ in front of the merger and deformation terms are small enough to absorb this prefactor, at least under our current choice of parameters.

One should also think of $\gamma \gg 1$, so that if $\|f\|_{\lambda, \gamma, \rho} \lesssim 1$, then $|f(s, B)| \lesssim \gamma^{-\area(s)}$, and we obtain the area decay this way. Thus, $\gamma^{-1}$ should be thought of as the rate of decay in the area law. 

The term $\rho^{B-K}$ is just needed to handle the ``base case" or i.e. the ``boundary term", where $s$ is a null loop. Ultimately, we want to compare $|\phi(\varnothing, K)|$ to our original partition function $|\phi(\varnothing, B)|$, and the additional weight $\rho^{B-K}$ allows us some room in making this comparison. The reason why we expect that it is safe to add this prefactor is because each deformation term comes with a tiny factor $\beta N^{-1}$, so that every plaquette we see, we pay a huge cost.
\end{remark}

The key property of the truncated model that will allow us to close the contraction mapping argument is the following result which bounds the total number of string operations. In particular, the total number of mergers is bounded.

\begin{lemma}\label{lemma:bound-on-loop-operations}
For any $(s, K)$ with $s$ non-null, any edge $\mbf{e}$ of $s$, and $K : \mc{P}_\Lambda^+ \ra \N$, we have that
\begin{equs}
| \mbb{D}_+(s, \mbf{e}, K) \cup \mbb{D}_-(s, \mbf{e}, K)| \leq 4d.
\end{equs}
If additionally $(s, K) \in \Omega_B$, then
\begin{equs}
| \mbb{S}_+(s, \mbf{e}) \cup \mbb{S}_-(s, \mbf{e})| &\leq 2dB, \\
| \mbb{M}_+(s, \mbf{e}) \cup \mbb{M}_-(s, \mbf{e})| &\leq 2dB. 
\end{equs}
\end{lemma}
\begin{remark}
In the large-$N$ limit, when the merger term in the master loop equation is not present, one can actually close a contraction estimate even without an upper bound on the number of splittings as in Lemma \ref{lemma:bound-on-loop-operations}, see \cite{Chatterjee2019a, BCSK2024}. This is done by using a mixed $\ell^1-\ell^\infty$ norm, where the $\ell^1$ norm comes from a sum over all ``loop profiles". Also, the deformation terms are always bounded by a dimensional constant. Thus, the main new thing here is that the number of merger terms is also bounded.
\end{remark}

\begin{proof}
Note that every edge $e \in \Z^d$ is contained in exactly $2(d-1)$ unoriented plaquettes. Thus,
\begin{equs}
| \mbb{D}_+(s, \mbf{e}, K) \cup \mbb{D}_-(s, \mbf{e}, K)|\leq 2 \times 2(d-1) < 4d.
\end{equs}
Next, let $(s,K) \in \Omega_{B}$ with $s=(s_1,\dots,s_n)$ non-null, and without loss of generality, let $\mbf{e}$ be an edge of $s_1$. Denote the corresponding lattice edge by $e \in E_\Lambda$. Let $k_i(e)$ denote the number of occurrences of $e$ or $e^{-1}$ in $s_i$. Then $\sum_{i=1}^{n} k_i(e) \leq n_e(s,K) \leq 2dB$, so
\begin{equs}
    | \mbb{S}_+(s, \mbf{e}) \cup \mbb{S}_-(s, \mbf{e})| &\leq k_1(e) \leq 2dB, \\
    | \mbb{M}_+(s, \mbf{e}) \cup \mbb{M}_-(s, \mbf{e})| &\leq \sum_{i=2}^{n}k_i(e) \leq 2dB,
\end{equs}
as desired.
\end{proof}

We are almost ready to carry out the contraction mapping argument. First, we define the operator appearing in the master loop equation.

\begin{definition}[Loop equation operator]
We define an operator $M$ on the space of functions $f : \Omega_B \ra \R$ as follows. Let $f : \Omega_B \ra \R$. Let $(s, K) \in \Omega_B$ be such that $s$ is non-null, and let $\mbf{e} = \mbf{e}_s$ be the first edge of $s$. Define
\begin{equs}
(M f)(s, K) := \mp \sum_{s' \in \mbb{S}_{\pm}(s, \mbf{e})} f(s', K) \mp \frac{1}{N^2} \sum_{s' \in \mbb{M}_{\pm}(s, \mbf{e})} f(s', K) \mp \frac{\beta}{N} \sum_{(s', K') \in \mbb{D}_{\pm}(s, \mbf{e}, K)} f(s', K').
\end{equs}
For null loops, we define
\begin{equs}
M f(\varnothing_x,K):= f(\varnothing_x ,K), ~~ x \in \Lambda.
\end{equs}
\end{definition}

By the master loop equation (Proposition \ref{prop:mle}), we have that $\phi = M \phi$, where $\phi$ is as defined in \eqref{eq:phi-s-K}.

\begin{prop}[Contraction estimate]\label{prop:contraction}
Let $\lambda, \rho \leq 1$, and $\gamma \geq 1$.
For any $f \in \R^{\Omega_B}$, we have that
\[ \|M f\|_{\lambda, \gamma, \rho} \leq \Big(2 d B \lambda + \frac{2d B}{\lambda N^2} + \frac{4d\beta \gamma}{\lambda \rho N} \Big) \|f\|_{\lambda, \gamma, \rho} + \sup_{x \in \Lambda} \sup_{K \leq B} \rho^{B-K} |f(\varnothing_x, K)|. \]
\end{prop}
\begin{proof}
Note that for any $x \in \Lambda$, the set $\{K : (\varnothing_x, K) \in \Omega_B\} = \{K : K \leq B\}$. Note also $\iota(\varnothing_x) = 0, \area(\varnothing_x) = 0$. Thus
\[ \|M f\|_{\lambda, \gamma, \rho} \leq \sup_{\substack{(s, K) \in \Omega_B \\ s \neq \varnothing}}  \lambda^{\iota(s)} \gamma^{\area(s)} \rho^{B-K} |f(s, K)| + \sup_{x \in \Lambda}\sup_{K \leq B} \rho^{B-K} |f(\varnothing_x, K)|.  \]
To finish, it suffices to show that
\[ \sup_{\substack{(s, K) \in \Omega_B \\ s \neq \varnothing}}  \lambda^{\iota(s)} \gamma^{\area(s)} \rho^{B-K} |f(s, K)| \leq \Big(2 d B \lambda + \frac{2d B}{\lambda N^2} + \frac{4d\beta \gamma}{\lambda \rho N} \Big) \|f\|_{\lambda, \gamma, \rho} .\]
Towards this end, fix $(s, K) \in \Omega_B$ with $s$ non-null. By the definition of $M$, it suffices to prove the following three estimates:
\begin{equs}
\lambda^{\iota(s)} \gamma^{\area(s)} \rho^{B-K} \sum_{s' \in \mbb{S}_{\pm}(s, \mbf{e}_s)} |f(s', K)| &\leq 2d B \lambda \|f\|_{\lambda, \gamma, \rho} \label{eq:intermediate-splitting}\\
\lambda^{\iota(s)} \gamma^{\area(s)} \rho^{B-K} \sum_{s' \in \mbb{M}_{\pm}(s, \mbf{e}_s)} |f(s', K)| &\leq \frac{2d B}{\lambda} \|f\|_{\lambda, \gamma, \rho} \label{eq:intermediate-merger}\\ 
\lambda^{\iota(s)} \gamma^{\area(s)} \rho^{B-K} \sum_{(s', K') \in \mbb{D}_{\pm}(s, \mbf{e}_s, K)} |f(s', K')| &\leq  \frac{4d \gamma}{\lambda \rho}  \|f\|_{\lambda, \gamma, \rho}. \label{eq:intermediate-deformation}
\end{equs}
\emph{Proof of \eqref{eq:intermediate-splitting}.}
Let $s' \in \mbb{S}_{\pm}(s, \mbf{e}_s)$. By Lemma \ref{lemma:splitting-complexity-after-splitting}, we have that $\iota(s') \leq \iota(s) - 1$, and so $\lambda^{\iota(s)} = \lambda \cdot \lambda^{\iota(s) - 1} \leq \lambda \cdot \lambda^{\iota(s')}$ (since $\lambda \leq 1$). We also have that $\area(s) = \area(s')$ (Lemma \ref{lemma:area-after-loop-operation}). Note also that $(s, K) \in \Omega_B$ implies that $(s', K) \in \Omega_B$. Thus 
\[\lambda^{\iota(s)} \gamma^{\area(s)} \rho^{B-K} |f(s', K)| \leq \lambda \lambda^{\iota(s')} \gamma^{\area(s')} \rho^{B-K} |f(s', K)| \leq \lambda \|f\|_{\lambda, \gamma, \rho}. \]
Since there are at most $2d B$ positive or negative splittings (Lemma \ref{lemma:bound-on-loop-operations}), the estimate \eqref{eq:intermediate-splitting} follows. \\

\noindent \emph{Proof of \eqref{eq:intermediate-merger}.} 
Let $s' \in \mbb{M}_{\pm}(s, \mbf{e}_s)$. By Lemma \ref{lemma:splitting-complexity-merger} and the fact that $\lambda \leq 1$, we have that $\lambda^{\iota(s)} = \lambda^{-1} \lambda^{\iota(s) + 1} \leq \lambda^{-1} \lambda^{\iota(s')}$. We also have that $\area(s) = \area(s')$ (Lemma \ref{lemma:area-after-loop-operation}). Note also that $(s, K) \in \Omega_B$ implies that $(s', K) \in \Omega_B$. Thus
\begin{equs}
\lambda^{\iota(s)} \gamma^{\area(s)} \rho^{B-K} |f(s', K)| \leq \lambda^{-1} \lambda^{\iota(s')} \gamma^{\area(s')} \rho^{B-K} |f(s', K)| \leq \lambda^{-1} \|f\|_{\lambda, \gamma, \rho}. 
\end{equs}
Since there are at most $2d B$ positive or negative mergers (Lemma \ref{lemma:bound-on-loop-operations}), the estimate \eqref{eq:intermediate-merger} follows. \\

\noindent \emph{Proof of \eqref{eq:intermediate-deformation}.}
Let $(s', K') \in \mbb{D}_{\pm}(s, \mbf{e}_s, K)$. By \eqref{eq:defomration-area-change}, we have that 
\[ \rho^{B-K} = \rho^{B - (K' + 1)} = \rho^{-1} \rho^{B-K'}. \]
By Lemma \ref{lemma:splitting-complexity-deformation}, we have that $\iota(s') \leq \iota(s) + 1$, and so $\lambda^{\iota(s)} \leq \lambda^{-1} \lambda^{\iota(s')}$ (as in the merger case). We also have that $\area(s) \leq \area(s') + 1$ (Lemma \ref{lemma:area-after-loop-operation}), and thus $\gamma^{\area(s)} \leq \gamma^{\area(s')+1}$ (since $\gamma \geq 1$). Finally, $(s, K) \in \Omega_B$ implies that $(s', K') \in \Omega_B$. Thus
\[ \lambda^{\iota(s)} \gamma^{\area(s)} \rho^{B-K} |f(s', K')| \leq \lambda^{-1} \rho^{-1} \gamma \lambda^{\iota(s')} \gamma^{\area(s')} \rho^{B-K'} |f(s', K')| \leq \frac{\gamma}{\lambda \rho} \|f\|_{\lambda, \gamma, \rho}.\]
Since there are at most $4d$ positive or negative deformations (Lemma \ref{lemma:bound-on-loop-operations}), we obtain the estimate \eqref{eq:intermediate-deformation}. This completes the proof of Proposition \ref{prop:contraction}.
\end{proof}

Before using the contraction estimate (Proposition \ref{prop:contraction}) to prove Theorem \ref{thm:truncated-model-area-law}, we first need to estimate the ``boundary term" by the partition function $Z_{\Lambda, \beta, N, B}$. In the following, recall that with $\phi$ defined as in \eqref{eq:phi-s-K}, the value $\phi(\varnothing_x, K)$ does not depend on the basepoint $x$, so that we abuse notation and write $\phi(\varnothing, K)$.

\begin{prop}\label{prop:partition-fn-bound}
With the choice of parameters as in \eqref{eq:N-large-assumption}-\eqref{eq:B-choice}, we have that
\begin{equs}
\sup_{K \leq B} e^{-(B-K)} |\phi(\varnothing, K)| \leq \phi(\varnothing, B) = Z_{\Lambda, \beta, N, B}. 
\end{equs}
\end{prop}
\begin{proof}
We previously already observed the equality (recall \eqref{eq:partition-fn-identity}), and thus we just show the inequality. We have that
\[ e^{-(B-K)} |\phi(\varnothing, K)| \leq \int dU \prod_{p \in \mc{P}_\Lambda^+} e^{-(B-K(p))} \big| \exp_{K(p)}(2 \beta \mrm{Re}\Tr(U_p))\big|. \]
By \eqref{eq:Truncation level comparison}, we have the pointwise bound
\begin{equs}
e^{-(B-K(p))} \big| \exp_{K(p)}(2 \beta \mrm{Re}\Tr(U_p))\big| \leq \exp_B(2 \beta \mrm{Re}\Tr(U_p)),
\end{equs}
and applying this for all $p \in \mc{P}_\Lambda^+$, we obtain the desired result.


\end{proof}

By combining the previous proposition with the contraction estimate (Proposition \ref{prop:contraction}), we can now show area law for the truncated model.

\begin{proof}[Proof of Theorem \ref{thm:truncated-model-area-law}] 
Let $\lambda = N^{-1}$, $\gamma = (10^3 d \beta)^{-1} \geq 1$, $\rho = e^{-1}$. Then by our assumptions on $N, \beta, B$ \eqref{eq:N-large-assumption}-\eqref{eq:B-choice}, we have by Propositions \ref{prop:mle} and \ref{prop:partition-fn-bound} that
\[ \|\phi\|_{\lambda, \gamma \rho} = \|M \phi\|_{\lambda, \gamma, \rho} \leq \frac{1}{2} \|\phi\|_{\lambda, \gamma, \rho} + Z_{\Lambda, \beta, N, B} .\]
In slightly more detail, we claim that
\begin{equs}
2d B \lambda + \frac{2d B}{\lambda N^2} + \frac{4d\beta \gamma}{\lambda \rho N} \leq \frac{1}{2}.
\end{equs}
This follows because $\lambda = N^{-1}$, $dB \ll N$ (which handles the first two terms), and
\begin{equs}
d \beta \gamma = 10^{-3} \ll 1, 
\end{equs}
which handles the third term. Letting $\ell$ be a rectangular loop, and recalling by Remark \ref{remark:configuration-space} that $(\ell, B) \in \Omega_B$, we obtain that
\[ \lambda^{\iota(\ell)} \gamma^{\area(\ell)} |\phi(\ell, B)| \leq  \|\phi\|_{\lambda, \gamma, \rho} \leq 2 \phi(\varnothing, B). \]
The desired result now follows upon recalling the values of $\lambda, \gamma$, using that $\iota(\ell) = |\ell|/4 - 1$, and noting that $\langle W_\ell \rangle_{\Lambda, \beta, N, B} = \phi(\ell, B) / Z_{\Lambda, \beta, N, B}$.
\end{proof}


\section{Original model}\label{section:original-model}

In this section, we prove Theorem \ref{thm:area-law-main} by adapting the contraction mapping argument of Section \ref{section:truncated-model}. 
The main result of this section is the following theorem. 
\begin{theorem}[Area Law]\label{thm:area-law}
Under the parameter assumptions on $N, \beta$ given in \eqref{eq:N-large-assumption}-\eqref{eq:beta-choice}, for any lattice $\Lambda$ and any rectangular loop $\ell$ in $\Lambda$, we have that
\begin{equs}
|\langle W_{\ell} \rangle_{\beta,\Lambda,N}| \leq 2  |\ell| N^{|\ell|/4-1} \alpha^{\area(\ell)},
\end{equs}
where 
\begin{equs}
\alpha = 2 \times 10^{3d} \max(10^3 d \beta, \exp(-10^{-7} d^{-2} N)) \ll 1.
\end{equs}
\end{theorem}

First, we show that Theorem \ref{thm:area-law} implies Theorem \ref{thm:area-law-main}. The main point is that in Theorem \ref{thm:area-law}, $N$ is assumed to be sufficiently large, whereas Theorem \ref{thm:area-law-main} is stated for all $N \geq 1$. To handle the latter case, we just use the fact that for each fixed $N$, we know area law at sufficiently small $\beta$.

\begin{proof}[Proof of Theorem \ref{thm:area-law-main}]
From \cite[Theorem 5.1]{osterwalderseiler1978}, for each $N \geq 1$, there exists some $\beta_0(d, N)$ such that for all $\beta \leq \beta_0(d, N)$, area law in the form stated in Theorem \ref{thm:area-law-main} holds. Combining this with Theorem \ref{thm:area-law}, we may define 
\begin{equs}
\beta_0(d) := \min\big(10^{-10d} d^{-1}, \min_{N \leq 10^{10d} d^{10}} \beta_0(d, N)\big),
\end{equs} which is then a threshold that works for all $N \geq 1$.
\end{proof}

We next discuss some elements in the proof of Theorem \ref{thm:area-law}, which is the main focus of the rest of Section \ref{section:original-model}. The main difficulty is that in principal, the total number of copies of a given plaquette may be unbounded, while the contraction mapping argument of Section \ref{section:truncated-model} only works if each plaquette has at most $\sim B$ copies. On the other hand, the point of choosing the parameters $\beta, B$ in the way specified by \eqref{eq:beta-choice}-\eqref{eq:B-choice} is that if a given plaquette has more than $B$ copies, then we are in the tail of the exponential, and this should already come with a huge cost (recall e.g. the estimate \eqref{eq:Truncated vs Untruncated exp}). 

\begin{figure}
\centering
\includegraphics[width=0.5\linewidth]{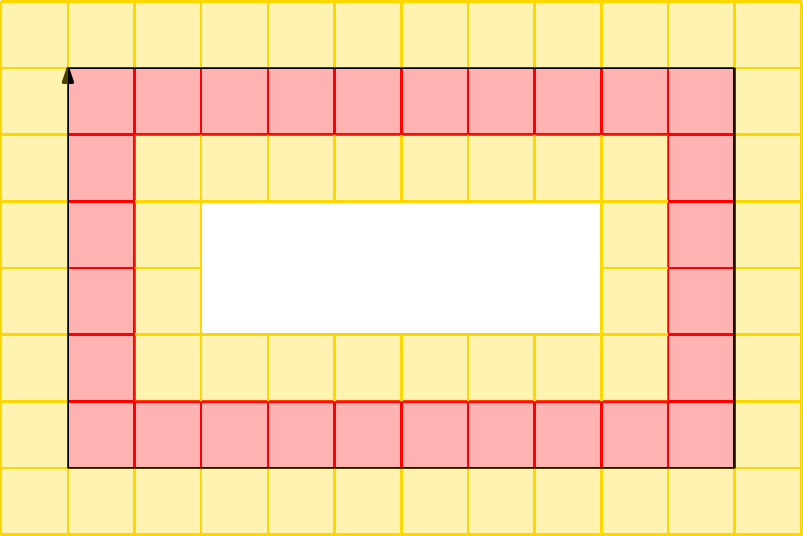}
\caption{Shown is our enemy scenario. It could happen that our loop $\ell$ (black) is completely contained in a set of bad plaquettes (red) which is of size $\sim |\ell|$. In this case, we only pay a cost which is exponential in $|\ell|$ due to the bad plaquettes, whereas we ultimately want to exhibit a cost which is exponential in $\area(\ell) \gg |\ell|$. The way we want to exhibit this cost is to begin exploring, and to say that we must see $\sim \area(\ell)$ many deformations, each of which comes with a cost. However, because every edge of the loop is contained in a bad plaquette, we cannot start our exploration. (We only want to explore at edges completely contained in good plaquettes, because only for these edges can we close the contraction mapping argument as in Section \ref{section:truncated-model}.) As described in the main text, the fix is to add another operation in the master loop equation, which introduces new loops corresponding to good plaquettes (yellow) which border the bad plaquette component in which we are currently stuck. We can then begin exploring at the edges of these new loops.}
\label{figure:loop-with-bad-plaquettes}
\end{figure}

However, even accounting for this, there is another key idea which is needed, which is dealing with the case where the loop is completely contained in the ``bad plaquettes", i.e. the plaquettes with too many copies. This is a problem because we only ever want to apply the master loop equation at edges which are only contained in good plaquettes (i.e., plaquettes with at most $B$ copies), so if every edge of the loop is contained in a bad plaquette, then there are no edges at which to apply the master loop equation. In terms of the surface exploration heuristic from Remark \ref{remark:surface-exploration}, we only want to explore at the edges of our string which are completely contained in good plaquettes, so if no such edges exist, then our exploration process is stuck. See Figure \ref{figure:loop-with-bad-plaquettes} for a visualization of such a situation. It may also be good to have the figure in mind while reading the ensuing discussion.

The way around this conundrum is to add another possible string operation to the master loop equation, which is applied precisely when the string is stuck in a set of bad plaquettes. The operation, which we term ``revival", simply proceeds by introducing new loops. These new loops are copies of good plaquettes that border the component of bad plaquettes in which the loop is currently stuck. The number of copies introduced is given by the current plaquette count, which recall tracks the total number of copies of the plaquette in our system. We can then continue exploration process by exploring the edges of these new loops. In terms of the surface exploration heuristic, the main point here is that we know there must be a surface spanning the original loop, and thus it must be possible to continue the exploration upon introducing the new loops, because at least one of the new loops must be part of the spanning surface. 

In the end, the area law decay will follow because in our exploration process, the combined number of bad plaquettes and good plaquettes that we see must be at least the area of the original loop (because there must be a surface spanning the loop). The bad plaquettes are extremely rare and thus they come with a huge cost, and the good plaquettes also come with a huge cost, essentially because in the master loop equation, deformations come with a prefactor of $\beta N^{-1}$.

We briefly summarize the rest of Section \ref{section:original-model}. In Section \ref{section:bad-plaquette-reduction}, we reduce Theorem \ref{thm:area-law} to proving a certain estimate for each fixed ``bad plaquette" configuration. In Section \ref{section:modified-mle}, we introduce a modified master loop equation, and in particular define the new revival operation. Finally, in Section \ref{section:contracton-mapping-original-model}, we prove an a priori estimate for solutions to the modified master loop equation, and then use this to prove Theorem \ref{thm:area-law}.

\subsection{Bad plaquette reduction}\label{section:bad-plaquette-reduction}

In this section, we perform a preliminary reduction towards the proof of Theorem \ref{thm:area-law}. The main result of this section is Proposition \ref{prop: reduction to fixed plaquette config}. Before we state the result, we include some preliminary discussion of plaquette sets.

\begin{convention}[Good and bad plaquettes]
In the following, we will often fix a plaquette set $P \sse \mc{P}_\Lambda^+$, which we will refer to as the set of ``bad plaquettes". Essentially, $P$ will track for us the plaquettes which have too many copies. Conversely, when we refer to a ``good plaquette", we mean a plaquette in $\mc{P}_\Lambda^+ \backslash P$, which will have at most $B$ copies. The good plaquettes are good because our contraction mapping argument from Section \ref{section:truncated-model} works at edges completely contained in good plaquettes.
\end{convention}

\begin{remark}
For a picture to keep in mind throughout the present section, it may help to think of the bad plaquette set $P$ as in Figure \ref{figure:loop-with-bad-plaquettes}. That is, it is a plaquette set which has size $\lesssim |\ell|$, and which may possibly contain many (or all) the edges of $\ell$ (which would prevent us from directly applying the argument of Section \ref{section:truncated-model}). This is our enemy scenario, because if, on the other hand, there are $\gtrsim \area(\ell)$ many bad plaquettes near our loop $\ell$, then we will automatically have the area law decay, because bad plaquettes are extremely rare.
\end{remark}

We will often talk about ``connected" plaquette sets, and so we define the following graph structure on the set of plaquettes.

\begin{definition}[Graph structure on plaquettes]
We say that two plaquettes $p,p' \in \mc{P}_{\Lambda}^+$ are neighbors if they share a common edge $e \in E_{\Lambda}$, and we write $p \sim p'$. Given a plaquette set $\mc{C} \sse \mc{P}_\Lambda^+$ and a loop $\ell$, we will say that $\mc{C} \cup \ell$ is connected if every connected component of $\mc{C}$ contains at least one edge of $\ell$. 
\end{definition}

A plaquette set $\mc{C}$ such that $\mc{C} \cup \ell$ is connected may contain several disjoint connected components. These components are connected to each other via the loop $\ell$.

\begin{definition}\label{def:plaquette-loop-cluster}
For any bad plaquette configuration $P \sse \mc{P}_{\Lambda}^+$, let 
\begin{equs}\label{eq:M-P}
M(P):= \max \big\{|P \cap \mc{C}|: \mc{C} \sse \mc{P}_\Lambda^+, ~\mc{C}\cup \ell \text{ is connected,}~ |\mc{C}|\leq \mathrm{area}(\ell) \big \}.
\end{equs}
\end{definition}

The following lemma bounds the number of sets $\mc{C}$ such that $\mc{C} \cup \ell$ is connected and which are at most a given size. Such statements are quite classical, but for completeness, we provide a proof. The reader may skip the proof on a first reading.

\begin{lemma}
Let $\ell$ be a rectangular loop. For any $M \geq 0$, we have that
\begin{equs}\label{eq:connected cluster bound}
\#\{\mc{C} \sse \mc{P}_\Lambda^+ :\mc{C}\cup \ell\text{ is connected}, |\mc{C}|\leq M\} \leq e^{2d|\ell|} 50^{dM}.
\end{equs}
\end{lemma}
\begin{proof}
The following is quite similar to the proof of \cite[Lemma 5.19]{BCSK2024}. It suffices to fix $1 \leq m \leq M$, and bound
\begin{equs}\label{eq:connected-bd-suffice}
\#\{\mc{C} \sse \mc{P}_\Lambda^+ :\mc{C}\cup \ell\text{ is connected}, |\mc{C}|= m\} \leq e^{2d|\ell|} 40^{dm}.
\end{equs}
The result then follows by summing in $1 \leq m \leq M$. For fixed $m$, any plaquette set $\mc{C}$ of size $|\mc{C}| = m$ such that $\mc{C} \cup \ell$ is connected decomposes into $k$ connected components, where $1 \leq k \leq m$. The sizes $n_1, \ldots, n_k \geq 1$ of these connected components form a partition of $m$, i.e. $n_1 + \cdots + n_k = m$. For $n \geq 1$, let $N(\ell, n)$ be the number of connected plaquette sets $\mc{C}$ which contain some edge of $\ell$. We then have that
\begin{equs}\label{eq:connected-plaquette-bound-intermediate}
\#\{\mc{C} \sse \mc{P}_\Lambda^+ :\mc{C}\cup \ell\text{ is connected}, |\mc{C}|= m\} \leq \sum_{k=1}^m \frac{1}{k!} \sum_{\substack{n_1, \ldots, n_k \geq 1 \\ n_1 + \cdots + n_k = m}} \prod_{i=1}^k N(\ell, n_i).
\end{equs}
Next, we claim that
\begin{equs}
N(\ell, n) \leq 2(d-1)|\ell| 10^{dn} \text{ for all $n \geq 1$}.
\end{equs}
This follows because each edge in $\ell$ is contained in at most $2(d-1)$ plaquettes, and for any fixed plaquette, the number of connected plaquette sets of size $n$ containing that plaquette is at most $10^{dn}$ (see e.g. \cite[Chapter 4]{GrimmettPercolation}, and in particular, (4.24) therein). Plugging this estimate into \eqref{eq:connected-plaquette-bound-intermediate} and combining with the elementary inequalities
\begin{equation}\label{eq:binomial-bound}
\dbinom{m-1}{k-1} \leq \bigg(\frac{(m-1) e}{k-1}\bigg)^{k - 1}, ~~\text{and}~~ \sup_{x > 0} (m/x)^x \leq e^{m/e},
\end{equation}
we obtain (here, note by stars-and-bars that the number of partitions $n_1 + \cdots + n_k = m$ such that $n_1, \ldots, n_k \geq 1$ is given by $\binom{m-1}{k-1}$)
\begin{equs}
\#\{\mc{C} \sse \mc{P}_\Lambda^+ :\mc{C}\cup \ell\text{ is connected}, |\mc{C}|= m\} &\leq \sum_{k=1}^m \frac{1}{k!} \binom{m-1}{k-1} (2(d-1)|\ell|)^k 10^{dm}  \\
&\leq \sum_{k=1}^m \frac{1}{k!} \bigg(\frac{(m-1) e}{k-1}\bigg)^{k - 1} (2(d-1)|\ell|)^k 10^{dm} \\
&\leq 10^{dm} e^{(m-1)/e} 
\sum_{k=1}^m \frac{1}{k!} (2e(d-1)|\ell|)^k \\
&\leq 10^{dm} e^{(m-1)/e} e^m e^{2(d-1)|\ell|}.
\end{equs}
The estimate \eqref{eq:connected-bd-suffice} now follows because $e^{1/e} \times e \times 10 < 40$.
\end{proof}

The following proposition is the main result of this section. It provides a sufficient condition under which area law holds. In the following, recall the truncated exponential and exponential tail functions $\exp_B, \tail_B$ from Definition \ref{def:truncated-exponential}. These were defined so that $\exp = \exp_B + \tail_B$.

\begin{prop}[Reduction to fixed bad plaquette configuration]\label{prop: reduction to fixed plaquette config}
Suppose that for some constants $C_1, C_2$, we have that for all $P \sse \mc{P}_\Lambda^+$,
\begin{equs}\label{eq:Fixed Plaquette Config Area Law}
\bigg|\int dU W_{\ell}(U) \prod_{p \in \mc{P}_{\Lambda}^+ \backslash P} &\exp_B(2 \beta \mathrm{Re}\mathrm{Tr}(U_p))\prod_{p \in P} \tail_B(2 \beta \mathrm{Re}\mathrm{Tr}(U_p)) \bigg| \\
& \leq C_1 \exp((B/40) M(P)) \exp(-C_2 (\area(\ell) - 20dM(P))_+)  \\
&\quad \quad \quad \quad ~\times\int dU \prod_{p \in \mc{P}_{\Lambda}^+ \backslash P} \exp_B(2 \beta \mathrm{Re}\mathrm{Tr}(U_p))\prod_{p \in P} \tail_B(2 \beta \mathrm{Re}\mathrm{Tr}(U_p)).
\end{equs}
Here, $(\area(\ell) - 20dM(P))_+ = \max(\area(\ell) - 20dM(P), 0)$. Then we have that
\begin{equs}\label{eq:deduce-area-law}
|\langle W_\ell \rangle_{\Lambda, \beta, N}| \leq C_1 e^{2d|\ell|} 50^{d \cdot \area(\ell)} 2^{\area(\ell)} \max(\exp(-C_2), \exp(-d^{-1} 10^{-3} B))^{\area(\ell)}.
\end{equs}
\end{prop}
\begin{remark}
Of course, the estimate in \eqref{eq:deduce-area-law} is only effective if $C_2 \gg d \log 50 + \log 2$. As we will see later, this will be guaranteed by our choice of parameters \eqref{eq:N-large-assumption}-\eqref{eq:B-choice}.
\end{remark}
\begin{proof}
First, for every plaquette configuration $P \sse \mc{P}_{\Lambda}^+$, let
\begin{equs}
\mc{C}(P) \in \{\mc{C}:\mc{C}\cup \ell \text{ is connected}, |\mc{C}|\leq \mathrm{area}(\ell)\}
\end{equs}
be a set which achieves maximal intersection, i.e. such that $|P \cap \mc{C}_P| = M(P)$. We may then estimate
\begin{equs}
\bigg|\int &dU W_{\ell}(U) \prod_{p \in \mc{P}_{\Lambda}^+} \exp(2 \beta \mathrm{Re} \mathrm{Tr}(U_p))\bigg|\\
&= \bigg|\int dU W_{\ell}(U) \prod_{p \in \mc{P}_{\Lambda}^+} \Big(\exp_B(2 \beta \mathrm{Re} \mathrm{Tr}(U_p))+ \tail_B(2 \beta \mathrm{Re} \mathrm{Tr}(U_p))\Big)\bigg|\\
&\leq \sum_{P \sse \mc{P}_{\Lambda}^+} \bigg|\int dU W_{\ell}(U) \prod_{p \in \mc{P}_{\Lambda}^+ \backslash P} \exp_B(2 \beta \mathrm{Re}\mathrm{Tr}(U_p))\prod_{p \in P} \tail_B(2 \beta \mathrm{Re}\mathrm{Tr}(U_p))\bigg|\\
&\leq \sum_{\mc{C}} \sum_{\substack{P \sse \mc{P}_\Lambda^+ \\ \mc{C}(P) = \mc{C}}} \bigg|\int dU W_{\ell}(U) \prod_{p \in \mc{P}_{\Lambda}^+ \backslash P} \exp_B(2 \beta \mathrm{Re}\mathrm{Tr}(U_p))\prod_{p \in P} \tail_B(2 \beta \mathrm{Re}\mathrm{Tr}(U_p))\bigg| \\
&\leq \sum_{\mc{C}} \sum_{\substack{P \sse \mc{P}_\Lambda^+ \\ \mc{C}(P) = \mc{C}}} \sum_{\substack{P_0 \sse \mc{C} \\ Q \sse \mc{P}_\Lambda^+ \backslash \mc{C} \\ P = P_0 \cup Q}}  \bigg|\int dU W_{\ell}(U) \prod_{p \in \mc{P}_{\Lambda}^+ \backslash P} \exp_B(2 \beta \mathrm{Re}\mathrm{Tr}(U_p))\prod_{p \in P} \tail_B(2 \beta \mathrm{Re}\mathrm{Tr}(U_p))\bigg|.
\end{equs}
Here, to be clear, the sum over $\mc{C}$ is a sum over $\mc{C} \sse \mc{P}_\Lambda^+$ such that $\mc{C} \cup \ell$ is connected and $|\mc{C}| \leq \area(\ell)$. Also, the innermost sum is really just a sum over a single pair $(P_0, Q)$, which is uniquely determined by the conditions $P_0 \sse \mc{C}, Q \sse \mc{P}_\Lambda^+ \backslash \mc{C}, P = P_0 \cup Q$. Also, note that since $\mc{C}(P) = \mc{C}$, we have that $|P_0| = M(P)$. Next, inserting the assumed bound \eqref{eq:Fixed Plaquette Config Area Law}, we obtain the further upper bound
\begin{equs}
\leq C_1 \sum_{\mc{C}} &\sum_{\substack{P \sse \mc{P}_\Lambda^+ \\ \mc{C}(P) = \mc{C}}} \sum_{\substack{P_0 \sse \mc{C} \\ Q \sse \mc{P}_\Lambda^+ \backslash \mc{C} \\ P = P_0 \cup Q}} \exp(-C_2(\area(\ell) - 20d |P_0|)_+) \exp((B/40)|P_0|) \\
&\times\int dU \prod_{p \in \mc{P}_{\Lambda}^+ \backslash P} \exp_B(2 \beta \mathrm{Re}\mathrm{Tr}(U_p))\prod_{p \in P} \tail_B(2 \beta \mathrm{Re}\mathrm{Tr}(U_p)).
\end{equs}
Next, by applying \eqref{eq:Truncated vs Untruncated exp} to bound $\tail_B \leq \exp(-B/10) \exp_B$ for each $p \in P_0$, we obtain the further upper bound
\begin{equs}
\leq C_1 \sum_{\mc{C}} &\sum_{\substack{P \sse \mc{P}_\Lambda^+ \\ \mc{C}(P) = \mc{C}}} \sum_{\substack{P_0 \sse \mc{C} \\ Q \sse \mc{P}_\Lambda^+ \backslash \mc{C} \\ P = P_0 \cup Q}} \exp(-C_2(\area(\ell) - 20d |P_0|)_+) \exp(-(B/20)|P_0|) \\
&\times\int dU \prod_{p \in \mc{P}_{\Lambda}^+ \backslash Q} \exp_B(2 \beta \mathrm{Re}\mathrm{Tr}(U_p))\prod_{p \in Q} \tail_B(2 \beta \mathrm{Re}\mathrm{Tr}(U_p)).
\end{equs}
Using that $(\area(\ell) - 20d |P_0|)_+ + 20d |P_0| \geq \area(\ell)$, and letting $\alpha := \max(\exp(-C_2), \exp(-d^{-1} 10^{-3} B))$, we may estimate
\begin{equs}
\exp(-C_2(\area(\ell) - 20d |P_0|)_+) \exp(-(B/20)|P_0|) \leq \alpha^{\area(\ell)}.
\end{equs}
Next, switching the order of summation so that the sum over $Q$ is outside, we further obtain
\begin{equs}
\leq C_1 \alpha^{\area(\ell)} \sum_{Q \sse \mc{P}_\Lambda^+} \sum_{\mc{C} \sse \mc{P}_\Lambda^+ \backslash Q} &\sum_{\substack{P \sse \mc{P}_\Lambda^+ \\ \mc{C}(P) = \mc{C}}} \sum_{\substack{P_0 \sse \mc{C} \\ P = P_0 \cup Q}} \int dU \prod_{p \in \mc{P}_{\Lambda}^+ \backslash Q} \exp_B(2 \beta \mathrm{Re}\mathrm{Tr}(U_p))\prod_{p \in Q} \tail_B(2 \beta \mathrm{Re}\mathrm{Tr}(U_p)).
\end{equs}
Now for any $\mc{Q}$ fixed, the number of possible choices of $\mc{C}$ is at most $e^{2d|\ell|} 50^{d\cdot \area(\ell)}$ by \eqref{eq:connected cluster bound}. Having fixed $\mc{C}$, the set $P$ is determined by the choice of $P_0 \sse \mc{C}$, and there are at most $2^{\area(\ell)}$ subsets of $\mc{C}$. It follows that the above is further upper bounded by
\begin{equs}
~&\leq C_1  e^{2d|\ell|} 50^{d\cdot \area(\ell)} 2^{\area(\ell)} \alpha^{\area(\ell)}\sum_{Q \sse \mc{P}_\Lambda^+} \int dU \prod_{p \in \mc{P}_{\Lambda}^+ \backslash Q} \exp_B(2 \beta \mathrm{Re}\mathrm{Tr}(U_p))\prod_{p \in Q} \tail_B(2 \beta \mathrm{Re}\mathrm{Tr}(U_p)) \\
&= C_1 e^{2d|\ell|} 50^{d\cdot \area(\ell)} 2^{\area(\ell)} \alpha^{\area(\ell)} Z_{\Lambda, \beta, N}.
\end{equs}
The desired result now follows by combining the previous chain of inequalities.
\end{proof}



\subsection{Modified master loop equation}\label{section:modified-mle}

In this section, we introduce a new string operation called ``revival", and then prove a modified master loop equation which involves this new operation. This will allow us to deal with the case where the current string is completely contained in bad plaquettes. The main results of this section are Propositions \ref{prop:mle1} and \ref{prop:mle2}, which when combined, give the modified master loop equation. Throughout this section, we fix a bad plaquette set $P \sse \mc{P}_\Lambda^+$, and set $M=M(P)$, where recall $M(P)$ is defined in \eqref{eq:M-P}. 



\begin{convention}[Good edges]
Having fixed the set of bad plaquettes $P$, we will say that an edge $e \in E_\Lambda$ is good if it is not contained in a bad plaquette, i.e. if it is completely contained in good plaquettes.
\end{convention}

\begin{definition}
For any subset of positively oriented plaquettes $Q$, let $o(Q)= \{p,p^{-1} : p \in Q\}$.
\end{definition}

\begin{definition}[Revival]
Let $s=(\ell_1,\dots,\ell_n)$ be a string, and $R$ be some collection of plaquettes $R \sse \mc{P}_{\Lambda}^+$ referred to as a \textit{revival set}. Let $J: o(R) \to \N$ be a plaquette count on $o(R)$. Suppose that $o(R)=\{p_1,\dots,p_m\}$, and that $J(p_i)=j_i$ for $1 \leq i \leq m$. We define the \textit{revival} as the string
\begin{equs}
    s \Delta_R J = (\ell^{(1)}_1,\dots,\ell_{j_1}^{(1)},\dots, \ell_1^{(m)},\dots,\ell_{j_m}^{(m)},\ell_1,\dots,\ell_n),
\end{equs}
where each $\ell_k^{(i)}$ is a copy of the plaquette loop $p_i$. The set of all revivals of $s$ with revival set $R$ will be denoted $\mc{R}(s,R)$.
\end{definition}

\begin{remark}
This new revival operation is the one that was mentioned in the discussion after the proof of Theorem \ref{thm:area-law-main} in the beginning of Section \ref{section:original-model}. As discussed there, the main point of this operation is to deal with the case where all edges of the string are contained in bad plaquettes (this is the notion of stuck string, which we will soon define in Definition \ref{def:stuck}). By introducing this new operation, we introduce new loops into the system, and we can continue our exploration from there (recall Figure \ref{figure:loop-with-bad-plaquettes}).
\end{remark}

As a matter of notation, we make the following definition involving operations on pairs $(s, K)$ of strings $s$ and plaquette counts $K$.

\begin{definition}[Operations on string and plaquette count pairs]\label{def:string-plaquette-operation}
Let $s \in \mc{S}_{\Lambda}$ be non-null and $K: \mc{P}_{\Lambda}^{+} \backslash P \to \N$ be a plaquette count supported on good plaquettes. We make the following definitions.
\begin{enumerate}
    \item If $s$ contains a good edge $\mbf{e}$, then $(s',K')$ is a positive merger of $(s,K)$ at $\mbf{e}$ if $s' \in \mbb{M}_+(s,\mbf{e})$ and $K'=K$. The set of all positive mergers of $(s,K)$ at $\mbf{e}$ is denoted by $\mbb{M}_+((s,K),\mbf{e})$.
    \item The sets $\mbb{M}_-((s,K),\mbf{e})$, $\mbb{S}_+((s,K),\mbf{e})$, and $\mbb{S}_-((s,K),\mbf{e})$ may be similarly defined.
    \item If $s$ contains a good edge $\mbf{e}$, we say that $(s',K')$ is a positive (resp. negative) deformation of $(s,K)$ at $\mbf{e}$, if we can express $s'=s \oplus_{\mbf{e}} p$  (resp. $s'=s \ominus_{\mbf{e}} p$) for some plaquette $p$, and moreover $K'$ 
    is given by (here $|p|$ denotes the positively oriented version of $p$, which may be $p$ itself if it is already positively oriented)
    \begin{equs}
    K'(q) = \begin{cases}
        K(q) & q \in \mc{P}_\Lambda^+ \backslash \{|p|\}, \\
        K(p) - 1 & q  = |p|.
        \end{cases}
    \end{equs}
    The set of all such $(s',K')$ is denoted $\mbb{D}_+((s,K),\mbf{e})$ (resp. $\mbb{D}_-((s,K),\mbf{e})$). 
    \item We say that $(s',K')$ is a revival of $(s,K)$ with revival set $R$ if $s' \in \mc{R}(s,R)$, $K'(p)=K(p)$ for $p \in \mc{P}_{\Lambda}^+ \backslash R$, $K'(p)=0$ for all $p \in R$, and the plaquette count $J : o(R) \ra \N$ such that $s' = s \Delta_R J$ satisfies $J(p) + J(p^{-1}) \leq K(p)$ for all $p \in \mc{P}_\Lambda^+$. In this case, we write $(s', K') = (s, K) \Delta_R J$. The set of all revivals of $(s,K)$ with revival set $R$ is denoted $\mc{R}((s,K),R)$.
\end{enumerate} 
\end{definition}

Next, we define what it means for a string to be ``stuck" in the bad plaquette set $P$.

\begin{definition}[Stuck strings]\label{def:stuck}
For a non-null string $s \in \mc{S}_{\Lambda}$, we say that $s$ is \textit{stuck} in $Q \sse P$ if $Q$ is the union of all connected components of $P$ touching $s$, and all edges in the loops of $s$ belong to plaquettes in $Q$.
\end{definition}

Next, we define the quantity for which we will derive the modified master loop equation. For any $s \in \mc{S}_{\Lambda}$ and plaquette count $K: \mc{P}_{\Lambda}^+ \backslash P \to \N$, let
\begin{equs}\label{eq:phi-s-K-original}
\phi(s,K):=\int dU &W_s(U) \prod_{p \in \mc{P}^+_{\Lambda}\backslash P}\exp_{K(p)}(2\beta \mathrm{Re} \mathrm{Tr}(U_p))\prod_{p \in P} \tail_B(2\beta \mrm{Re}\Tr(U_p)).
\end{equs}
Next, we state the first main result of this section.

\begin{prop}[Master loop equation when not stuck]\label{prop:mle1}
Let $s \in \mc{S}_\Lambda$ be non-null with at least one good edge $\mbf{e}$, and let $K : \mc{P}_\Lambda^+ \ra \N$ be a plaquette count. We have that
\begin{equs}
\phi(s, K) = \mp \sum_{s' \in \mbb{S}_{\pm}(s, \mbf{e})} \phi(s', K) \mp \frac{1}{N^2} \sum_{s' \in \mbb{M}_{\pm}(s, \mbf{e})} \phi(s', K) \mp \frac{\beta}{N} \sum_{(s', K') \in \mbb{D}_{\pm}(s, \mbf{e}, K)} \phi(s', K').
\end{equs}
\end{prop}
\begin{proof}
We apply Proposition \ref{prop:general-mle} with the choices  $\rho_p(x)=\exp_{K(p)}(x)$ for $p \in \mc{P}_{\Lambda}^+ \backslash P$, and $\rho_p(x)=\tail_B$ for $p \in P$, and use the fact that $\frac{d}{dx}\exp_k(x)=\exp_{k-1}(x)$.
\end{proof}

Before we state the second main result of this section (the master loop equation for stuck strings), we make the following definitions.

\begin{definition}[Plaquette boundary]
For a set of bad plaquettes $Q \sse P$, we define the boundary
\begin{equs}
    \partial Q:=\{p \in \mc{P}_{\Lambda}^+ \backslash P: p \sim p' \text{ for some } p' \in Q\},
\end{equs}
i.e. $\ptl Q$ is the set of good plaquettes which share an edge with an element of $Q$.
\end{definition}

Next, we define a modified notion of area for strings, which will be slightly more convenient in the following. We remark that for rectangular loops, this modified notion of area coincides with our previous definition of area (Definition \ref{def:area}), which coincides with the usual notion of the area of a rectangle.

\begin{definition}[Modified area]\label{def:modarea}
For a string $s \in \mc{S}_{\Lambda}$, define the modified area of $s$ as 
\begin{equs}
\modarea(s):= \min \{|\mathrm{supp}(J)|: J: \mc{P}_{\Lambda} \to \N \text{ such that } (s,J) \text{ is balanced}\}.
\end{equs}
\end{definition}

In the rest of this section, whenever we refer to area, by default we will mean $\modarea$. We remark that the reason why $\modarea(\ell) = \area(\ell)$ for a rectangular loop $\ell$ is because the plaquette set $J$ which makes $(\ell, J)$ balanced and which minimizes $\sum_p J(p)$ is $0/1$-valued.


\begin{prop}[Master loop equation for stuck strings]\label{prop:mle2}
Let $s \in \mc{S}_\Lambda$ be a non-null string which is stuck in $Q \sse P$. Let $K : \mc{P}_\Lambda^+ \backslash P \ra \N$ be a plaquette count. Suppose that $|Q| < \modarea(s)$. Then
\begin{equs}
\phi(s,K)=\sum_{\substack{ J: \mc{P}_{\Lambda} \to \N \\ \varnothing \neq \mathrm{supp}(J) \sse \partial Q \\ J \leq K}} \frac{(N\beta)^J}{J!} \phi((s,K) \Delta_{\partial Q} J),
\end{equs}
where by $J \leq K$ we mean that for any oriented plaquette $p$, we have that $J(p)+J(p^{-1}) \leq K(|p|)$ (here $|p|$ denotes the positively oriented version of $p$). Note the RHS above is $0$ if $K\equiv 0$ on $\partial Q$.
\end{prop}

Before proving Proposition \ref{prop:mle2}, we collect a few key properties of balanced pairs $(s,K)$ (recall the definition of balanced in Definition \ref{def:balanced}).

\begin{lemma}\label{lemma:not-balanced-zero}
If $(s,J)$ is not balanced, then
\begin{equs}
\int dU W_s(U) \prod_{p \in \mc{P}_{\Lambda}} \mathrm{Tr}(U_p)^{J(p)} = 0.
\end{equs}
\end{lemma}
\begin{proof}
By assumption, there is an edge $e_0 \in E_\Lambda^+$ such that $m_{e_0} := n_{e_0}(s, J) - n_{e_0^{-1}}(s, J) \neq 0$, where recall $n_{e_0}(s, J)$ is the total number of appearances of $e_0$ in $(s, J)$, and similarly $n_{e_0^{-1}}(s, J)$ is the total number of appearances of $e_0^{-1}$ in $(s, J)$. Now take any $\theta \in [0, 2\pi)$ such that $e^{\icomplex \theta m_{e_0}} \neq 1$. By Haar invariance, $e^{\icomplex \theta} U_{e_0}$ is still Haar distributed, and thus the lattice gauge configuration $\tilde{U}$ defined by
\begin{equs}
\tilde{U}_e :=  \begin{cases} U_e & e \neq e_0 \\ e^{\icomplex \theta} U_{e_0} & e = e_0 
\end{cases}
\end{equs}
is still product Haar distributed, and thus we obtain
\begin{equs}
\int dU W_s(U) \prod_{p \in \mc{P}_{\Lambda}} \mathrm{Tr}(U_p)^{J(p)} = e^{\icomplex \theta m_e} \int dU W_s(U) \prod_{p \in \mc{P}_{\Lambda}} \mathrm{Tr}(U_p)^{J(p)}.
\end{equs}
Since $e^{\icomplex \theta m_e} \neq 1$, we obtain that the integral must be zero, as desired.
\end{proof}

\begin{lemma}
Suppose that $(s,J_1+J_2)$ and $(\varnothing_x,J_2)$ are both balanced (the second condition is independent of the choice of $x \in \Lambda$). Then $(s,J_1)$ is balanced.
\end{lemma}

\begin{proof}
For any edge $e \in E_{\Lambda}^+$, $e$ and $e^{-1}$ appear the same number of times in $s$ and $J_1+J_2$ together, as well as just in $J_2$. Thus $e$ and $e^{-1}$ appear the same number of times in $s$ and $J_1$ together.
\end{proof}

\begin{lemma}\label{lemma:balanced-support-at-least-area}
Suppose that $(s,J)$ is balanced, and let $S$ be the union of all connected components of $\mathrm{supp}(J)$ which neighbor any loop of the string $s$. Then $|S| \geq \modarea(s)$.
\end{lemma}

\begin{proof}
We can decompose $J=J_1+J_2$ where $J_1:=J \mathbbm{1}_{S}$. Since $(s, J_1+J_2)$ is balanced and any edge intersecting the support of $J_2$ does not lie on any loop of $s$, $(\varnothing,J_2)$ is balanced. Therefore $(s,J_1)$ is balanced, and thus necessarily $S = \mathrm{supp}(J_1)$ satisfies $|S| \geq \modarea(s)$.
\end{proof}

\begin{proof}[Proof of Proposition \ref{prop:mle2}]
We begin by Taylor expanding the truncated exponentials at every $p \in \partial Q$, so that

\begin{equs}
    \phi(s,K) &= \sum_{\substack{J: \mc{P}_{\Lambda} \to \N \\ \mathrm{supp}(J) \sse \partial Q \\ J \leq K}} \int dU W_s(U)  \prod_{p \in o(\partial Q)}\frac{\beta^{J(p)} \mathrm{Tr}(U_p)^{J(p)}}{J(p)!}\prod_{p \in \mc{P}_{\Lambda} \backslash P}\exp_{K'(p)}(2 \beta \mrm{Re}\Tr(U_p))\prod_{p \in  P}\tail_B(2 \beta \mrm{Re}\Tr(U_p))\\
    &=\sum_{\substack{ J: P_{\Lambda} \to \N \\ \mathrm{supp}(J) \sse \partial Q \\ J \leq K}} \frac{(N\beta)^J}{J!} \phi((s,K) \Delta_{\partial Q}J),
\end{equs}
where in the right hand side of the first line, $K'(p) =0$ for all $p \in \partial Q$, and $K'(p)=K(p)$ otherwise. To conclude, we just need to prove that the term above corresponding to $J\equiv0$ is $0$. To see this, we expand the remaining exponentials in series so that,
\begin{equs}
\int dU W_s(U) \prod_{p \in \mc{P}_{\Lambda} \backslash P} &\exp_{K'(p)}(2 \beta \mrm{Re}\Tr(U_p))\prod_{p \in  P}\tail_B(2 \beta \mrm{Re}\Tr(U_p)) \\
&=
\sum_{\Tilde{K}} \frac{\beta^{\Tilde{K}}}{\Tilde{K}!}\int dU W_s(U) \prod_{p \in \mc{P}_{\Lambda}} \mathrm{Tr}(U_p)^{\Tilde{K}(p)} \label{eq:J=0 term},
\end{equs}
where the $\Tilde{K}$ in the sum above are $0$ on $\partial Q$, $\Tilde{K} \leq K$ on the remaining good plaquettes, and $\Tilde{K}>B$ on the bad plaquette set $P$. The key observation is that for any such $\Tilde{K}$, $(s,\Tilde{K})$ cannot be balanced. To see this, let  $S$ be the union of all connected components of $\mathrm{supp}(\Tilde{K})$ neighboring $s$. Then we have that $S \sse Q$, because $\tilde{K} = 0$ on $\ptl Q$. By Lemma \ref{lemma:balanced-support-at-least-area} combined with the fact that $|Q| <\modarea(s)$, we have that the pair $(s, \Tilde{K})$  cannot be balanced. As a result, every integral in the sum \eqref{eq:J=0 term} must vanish (by Lemma \ref{lemma:not-balanced-zero}), thereby concluding the proof.
\end{proof}

\subsection{A priori estimate}\label{section:contracton-mapping-original-model}

In this section, we prove an a priori estimate for solutions to the modified master loop equation from Section \ref{section:modified-mle}, which will then allow us to obtain the estimate assumed in Proposition \ref{prop: reduction to fixed plaquette config}, which then will allow us to prove Theorem \ref{thm:area-law}. The main proposition in the section is Proposition \ref{prop:fixed-point-estimate}, which is the aforementioned a priori estimate, which we will gradually build towards.

Throughout this section, fix a rectangular loop $\ell$ belonging to a lattice $\Lambda$. We also fix a subset $P \sse \mc{P}_\Lambda^+$ of bad plaquettes. We will consider triples $(s, K, Q)$ where $s \in \mc{S}_{\Lambda}$, $K: \mc{P}_{\Lambda}^+ \backslash P \to \N$, and $Q \sse P$.

\begin{definition}[String operations on triples]
Suppose $s \in \mc{S}_{\Lambda}$, $K: \mc{P}_{\Lambda}^{+} \backslash P \to \N$, and $Q \sse P$. We make the following definitions. 
\begin{enumerate}
\item If $s$ contains a good edge $\mbf{e}$, then $(s',K',Q')$ is a positive merger of $(s,K,Q)$ at $\mbf{e}$ if $(s',K') \in \mbb{M}_+((s,K),\mbf{e})$ and $Q'=Q$. The set of all positive mergers of $(s,K)$ at $\mbf{e}$ is denoted by $\mbb{M}_+((s,K,Q),\mbf{e})$. We can define $\mbb{M}_-((s,K,Q),\mbf{e})$, $\mbb{S}_{\pm}((s,K,Q),\mbf{e})$, and $\mbb{D}_{\pm}((s,K,Q),\mbf{e})$ similarly. 
\item If $s$ is stuck in a minimal collection of connected components $\tilde{Q} \sse P$, then $(s',K',Q')$ is a revival of $(s,K,Q)$ if (1) $(s',K') \in \mc{R}((s,K),\partial \tilde{Q})$, (2) the plaquette count $J$ such that $(s', K') = (s, K) \Delta_{\ptl Q} J$ satisfies $J \neq 0$, and (3) $Q'=\tilde{Q} \cup Q$. The set of all revivals of $(s,K,Q)$ is denoted $\mc{R}(s,K,Q)$.
\end{enumerate}
\end{definition}

\begin{remark}\label{remark:edge-total-monotone}
Note that for any $(s', K', Q')$ obtained by applying a string operation to $(s, K, Q)$, we have that $n_e(s', K') \leq n_e(s, K)$ for all good edges $e$. In words, the total number of copies of any given good edge is non-increasing under the application of string operations.

Additionally, from the definition of revivals (Definition \ref{def:string-plaquette-operation}), we have that for all $(s', K', Q') \in \mc{R}(s, K, Q)$, the plaquette count $K'$ is zero on $\ptl (Q' \backslash Q)$. In words, this is because to obtain $(s', K', Q')$ from $(s, K, Q)$, we introduce a number of loops corresponding to plaquettes in $\ptl(Q' \backslash Q)$, while simultaneously zeroing out the plaquette count on $\ptl (Q' \backslash Q)$.
\end{remark}





\begin{convention}
In this section, our plaquette count $K : \mc{P}_\Lambda^+ \backslash P \ra \N$ will always be supported on good plaquettes, even if we do not explicitly say so. Thus, when we write $B - K$, we mean the plaquette count supported on $\mc{P}_\Lambda^+ \backslash P$ whose value at any given good plaquette $p$ is given by $B - K(p)$. In particular, the value of $B - K$ at a bad plaquette is zero. In particular, when we use $B$ to denote a plaquette count, we mean the plaquette count which takes the constant value $B$ on the good plaquettes, and is zero on the bad plaquettes.
\end{convention}

Next, we show how the area (from Definition \ref{def:modarea}) changes under a string operation.

\begin{lemma}\label{lemma:area-after-loop-operation-original}
Let $s \in \mc{S}_\Lambda$ be non-null, and $\mbf{e}$ be an edge of $s$. For all $s' \in \mbb{S}_{\pm}(s, \mbf{e}) \cup \mbb{M}_{\pm}(s, \mbf{e})$, we have that $\modarea(s') = \modarea(s)$. For $s' \in \mbb{D}_{\pm}(s, \mbf{e})$, we have that $\modarea(s) \leq \modarea(s') + 1$. Lastly, if $(s',K',Q') \in \mc{R}(s,K,Q)$, then $\modarea(s) \leq \modarea(s')+8(d-1)(|Q'|-|Q|)$.
\end{lemma}
\begin{proof}
From the definition of $\modarea$ as well as the definitions of splittings and mergers, we have that for any string $s' \in \mbb{S}_{\pm}(s, \mbf{e}) \cup \mbb{M}_{\pm}(s, \mbf{e})$ and any plaquette count $J$, the pair $(s, J)$ is balanced if and only if $(s', J)$ is balanced. This implies that $\modarea(s') = \modarea(s)$ for $s'$.

In the case of deformations, suppose $s'$ is a deformation of $s$ using a plaquette $p$. Then for any plaquette count $J'$ such that $(s', J')$ is balanced, we have that $(s, J' + \delta_{p^{-1}})$ is balanced and $|\mathrm{supp}(J' + \delta_{p^{-1}})| \leq |\mathrm{supp}(J')|+1 $. Taking $J'$ to be minimal area, we obtain the inequality $\modarea(s) \leq \modarea(s') + 1$. 

Finally, suppose that $(s',K',Q') \in \mc{R}(s,K,Q)$ and $s'=s \Delta_{\partial (Q' \backslash Q)} J$ for some plaquette count $J$. For each plaquette, there are at most $4 \times 2(d-1)=8(d-1)$ neighboring plaquettes, and thus we can crudely bound the support of $J$ by $|\partial (Q' \backslash Q)| \leq 8(d-1)|Q' \backslash Q| = 8(d-1)(|Q'|-|Q|)$. Moreover, if $(s',\tilde{J})$ is balanced, then $(s,\tilde{J}+J)$ is balanced, and $|\mathrm{supp}(\tilde{J}+J)| \leq |\mathrm{supp}(\tilde{J})|+|\mathrm{supp}(J)| \leq |\mathrm{supp}(\tilde{J})|+8(d-1)(|Q'|-|Q|)$. Taking $\tilde{J}$ such that $|\mathrm{supp}(\tilde{J})|=\modarea(s')$, we obtain the desired result.
\end{proof}

Next, we begin to define the configuration space $\Omega$ which will serve as the analog of $\Omega_B$ from Definition \ref{def:truncated-space}. First, we define the notion of string trajectory. In the following, for every non-null string $s \in \mc{S}_{\Lambda}$ which is not stuck in $P$, we fix a particular good edge $\mbf{e}=\mbf{e}_{s}$ of $s$. 

\begin{definition}[String trajectory]
A \textit{string trajectory} starting at $(s,K,Q)$ of length $l$ is a finite sequence
\begin{equs}
\Big((s_0,K_0,Q_0), (s_1,K_1,Q_1), \ldots, (s_l,K_l,Q_l)\Big),
\end{equs}
where $(s_0,K_0,Q_0)=(s,K,Q)$, and at each non-terminal step $i<l$, $s_i$ is non-null, and if $s_i$ is stuck in a minimal union of connected components $Q \sse P$, then $|Q|<\modarea(s_i)$ and
\begin{equs}
(s_{i+1},K_{i+1},Q_{i+1}) \in \mc{R}(s_i,K_i,Q_i).
\end{equs}
In the case where $s_i$ is not stuck, then
\begin{equs}
(s_{i+1},K_{i+1},Q_{i+1}) \in \mbb{S}_{\pm}((s_i,K_i,Q_i),\mbf{e})\cup \mbb{M}_{\pm}((s_i,K_i,Q_i),\mbf{e})\cup\mbb{D}_{\pm}((s_i,K_i,Q_i),\mbf{e}).
\end{equs}
\end{definition}
One should think of a string trajectory as encoding an ``exploration process", where at each step of this process, we either choose a good edge of the string to explore, or if the string is stuck, we introduce new plaquettes through the revival operation.

Next, we prove the following lemma which says that in this exploration process, we are always tracing out connected components. In the following, it may help to recall the definition \eqref{eq:M-P} of $M(P)$.

\begin{lemma}[String trajectories trace out cluster]\label{lemma:string-trajectory-cluster}
For any string trajectory 
\begin{equs}
\Big((s_0,K_0,Q_0), (s_1,K_1,Q_1), \ldots ,(s_l,K_l,Q_l)\Big)
\end{equs}
starting at $(\ell,B,\varnothing)$, for each $0 \leq i \leq l$, we have that 
\begin{equs}\label{eq:trace-out-cluster}
\mathrm{supp}(B-K_i)\cup Q_i \in \{\mc{C} \sse \mc{P}_{\Lambda}^+:  \mc{C} \cup \ell \text{ is connected}\}.
\end{equs}
If $s_i$ is stuck in a minimal union of connected components $\tilde{Q}_i$, we also have that 
\begin{equs}\label{eq:trace-out-cluster-2}
\mathrm{supp}(B-K_i)\cup Q_i \cup \tilde{Q}_i \in \{\mc{C} \sse \mc{P}_{\Lambda}^+: \mc{C} \cup \ell \text{ is connected}\}.
\end{equs}
We call elements of $\{\mc{C} \sse \mc{P}_{\Lambda}^+: \mc{C} \cup \ell \text{ is connected}\}$ connected clusters of $\ell$.
\end{lemma}
\begin{proof}
We first show the first claim, i.e. \eqref{eq:trace-out-cluster}. We proceed by induction. We claim that for all $0 \leq i \leq l$, we have that \eqref{eq:trace-out-cluster} is true, and also that
\begin{equs}\label{eq:string-in-cluster}
s_i \sse  \mrm{supp}(B-K_i) \cup Q_i \cup \ell.
\end{equs}
When $i = 0$, by assumption, $(s_0, K_0, Q_0) = (\ell, B, \varnothing)$, and $\mrm{supp}(B - B) = \varnothing$, so that $\mrm{supp}(B - B) \cup Q_0 = \varnothing \in \{\mc{C} : \mc{P}_\Lambda^+ : \mc{C} \cup \ell \text{ is connected}\}$. Also, $s_0 = \ell \sse \ell$, and thus \eqref{eq:string-in-cluster} is also true.

Next, suppose that the first claim is true for all $1 \leq i \leq j$, for some $1 \leq j < l$. If $s_{j+1}$ is a splitting or merger of $s_j$, then $K_{j+1} = K_j$ and $Q_{j+1} = Q_j$, and thus $\mrm{supp}(B - K_{j+1}) \cup Q_{j+1} = \mrm{supp}(B - K_j) \cup Q_j$, which shows \eqref{eq:trace-out-cluster} in this case. Moreover, $s_{j+1} \sse s_j$ (in terms of unoriented edge sets), and thus \eqref{eq:string-in-cluster} is also true in this case.

If $s_{j+1}$ is a deformation of $s_j$, then $\mrm{supp}(B-K_{j+1})$ at most may only differ from $\mrm{supp}(B-K_j)$ at a single plaquette $p_0 \in \mc{P}_\Lambda^+$, and moreover, $p_0$ must contain an edge of $s_j$. Also, $Q_{j+1} = Q_j$. Let $\mc{C}_j \sse \mc{P}_\Lambda^+$ be such that $\mc{C}_j \cup \ell$ is connected, and such that
\begin{equs}
\mrm{supp}(B-K_j) \cup Q_j = \mc{C}_j.
\end{equs}
Since \eqref{eq:string-in-cluster} holds for $s_j$, we then have that $\mc{C}_i \cup \{p_0\} \cup \ell$ is connected, and that
\begin{equs}
\mrm{supp}(B-K_{j+1}) \cup Q_{j+1} = \mc{C}_i \cup \{p_0\}.
\end{equs}
Thus, \eqref{eq:trace-out-cluster} holds. Finally, from the definition of deformation, we also have that 
\begin{equs}
s_{j+1} \sse s_j \cup p_0 \sse \mrm{supp}(B-K_{j+1}) \cup Q_{j+1}, 
\end{equs}
which shows \eqref{eq:string-in-cluster}.

Finally, suppose that $s_{j+1}$ is a revival of $s_j$. Then $s_j$ is stuck in a minimal union of connected components $\tilde{Q}_j \sse P$, and $Q_{j+1} = Q_j \cup \tilde{Q}_j$.  Also, there is some $R \sse \ptl \tilde{Q}_j$ such that $\mrm{supp}(B-K_{j+1}) = \mrm{supp}(B-K_j) \cup R$, and $s_{j+1} \sse s_j \cup R$. Since \eqref{eq:string-in-cluster} holds for $s_j$, we then obtain
\begin{equs}
s_{j+1} \sse s_j \cup R \sse \mrm{supp}(B-K_{j}) \cup R \cup Q_i \sse \mrm{supp}(B -K_{j+1}) \cup Q_{i+1},
\end{equs}
and thus \eqref{eq:string-in-cluster} holds for $s_{j+1}$. Finally, let $\mc{C}_j \sse \mc{P}_\Lambda^+$ be such that $\mc{C}_j \cup \ell$ is connected, and such that
\begin{equs}
\mrm{supp}(B-K_j) \cup Q_j = \mc{C}_j.
\end{equs}
Note that each component of $R$ is connected to $\tilde{Q}_j$, and also each component of $\tilde{Q}_j$ is connected to some edge of $s_j$, and thus since $s_j \sse \mc{C}_j \cup \ell$, it follows that $\mc{C}_j \cup R \cup \tilde{Q}_j \cup \ell$ is connected, and
\begin{equs}
\mrm{supp}(B-K_{j+1}) \cup Q_{j+1} = \mrm{supp}(B-K_j) \cup Q_j \cup R \cup \tilde{Q}_j = \mc{C}_j \cup R \cup \tilde{Q}_j,
\end{equs}
and thus \eqref{eq:trace-out-cluster} holds. This concludes the proof of the inductive step, and thus we have shown \eqref{eq:trace-out-cluster}.

To see \eqref{eq:trace-out-cluster-2}, let $\mc{C}_i \sse \mc{P}_\Lambda^+$ be such that $\mc{C}_i \cup \ell$ is connected, and such that
\begin{equs}
\mrm{supp}(B-K_i) \cup Q_i = \mc{C}_i.
\end{equs}
Then $s_i \sse \mc{C}_i$ (by \eqref{eq:string-in-cluster}), and by definition, each component of $\tilde{Q}_i$ is connected to an edge of $s_i$. Consequently, $\mc{C}_i \cup \tilde{Q}_i \cup \ell$ is connected, which is the claim \eqref{eq:trace-out-cluster-2}.
\end{proof}

\begin{definition}[Configuration space]\label{def:configuration-space}
Let $\Omega$ be the set of all triples $(s,K,Q)$ which appear in a string trajectory starting at $(\ell,B,\varnothing)$, and such that $|\mathrm{supp}(B-K) \cup Q| \leq \modarea(\ell)$. Here, technically $\Omega = \Omega_{\ell, B, P, \Lambda}$ depends on $\ell, B, P, \Lambda$, but as these are all fixed for us, we will omit the dependence of $\Omega$ on these parameters. Note that $\Omega$ is finite.
\end{definition}

\begin{remark}\label{remark:bound-on-loop-operations-string-trajectory}
Recalling that the total number of copies of any given good edge is non-increasing under application of string operations (Remark \ref{remark:edge-total-monotone}), we have that for all $(s, K, Q) \in \Omega$, for any good edge $e$, the number of copies $n_e(s, K) \leq n_e(\ell, B) \leq 2dB$. From this, it follows that the total number of mergers and splittings of $(s, K, Q)$ is at most $2dB$ (as in Lemma \ref{lemma:bound-on-loop-operations}).
\end{remark}

Next, we define the notion of terminal string, which should be thought of as the boundary of the configuration space $\Omega$. In particular, if our exploration process has reached a terminal string, we stop exploring. Remark \ref{remark:stop-exploring} gives some intuition as to why it suffices to stop exploring at terminal strings.

\begin{definition}[Terminal strings]
We call the triple $(s,K,Q) \in \Omega$ \textit{terminal} if any string trajectory $\big((s_0,K_0,Q_0), (s_1,K_1,Q_1), \ldots, (s_l,K_l,Q_l)\big)$ starting at $(s,K,Q)$ necessarily has $(s_1,K_1,Q_1) \notin \Omega$, or if there are simply no such string trajectories of length greater than $0$. The set of all terminal $(s,K,Q) \in \Omega$ is denoted $\partial\Omega$, and we set $\Omega^\circ :=\Omega \backslash \partial \Omega$.
\end{definition}

\begin{remark}[Examples of terminal strings]
We remark that any triple of the form $(\varnothing_x, K, Q)$ is a terminal string, because there are no string operations on null loops. Another example of a terminal string includes a triple $(s, K, Q)$ where $s$ is stuck in a minimal union of connected components $\tilde{Q} \sse Q$, such that $|\mrm{supp}(B-K) \cup Q \cup \tilde{Q}| > \modarea(\ell)$. 
\end{remark}

We next begin towards Proposition \ref{prop: terminal string decomposition}, which gives a decomposition of the set of terminal strings. Before we state and prove the result, we need some preliminary discussion.


\begin{definition}
Let $s$ be a string and $K$ be a plaquette count.
We say that the pair $(s,K)$ has an \textit{isolated edge} $e \in E^+(\Lambda)$, if $n_e(s, K) = 1$, i.e. if the number of occurrences of $e$ and $e^{-1}$ in $(s, K)$ is exactly 1.
\end{definition}

Note that by definition, any pair $(s, K)$ which contains an isolated edge cannot be balanced. 



\begin{definition}
Let $\ptl_0 \Omega \sse \ptl \Omega$ be the set of $(s, K, Q) \in \ptl \Omega$ such that either (1) $(s, K)$ contains an isolated edge, or (2) $s$ is stuck in some minimal union of connected components $\tilde{Q} \sse \mc{P}_\Lambda^+$, $|\tilde{Q}| < \modarea(s)$, and $\mc{R}(s, K, Q) = \varnothing$. Let $\ptl_1 \Omega := \ptl \Omega \backslash \ptl_0 \Omega$, so that $\ptl \Omega = \ptl_0 \Omega \sqcup \ptl_1 \Omega$.
\end{definition}


\begin{remark}\label{remark:ptl-0-omega}
In terms of the surface exploration heuristic, the set $\ptl_0 \Omega$ corresponds to the situation where in the middle of the surface exploration, we realize that it is not possible to see a surface spanning the original rectangular loop $\ell$. We should then be able to conclude that this case contributes zero. This is part Proposition \ref{prop: terminal string decomposition}, to be stated shortly. See Remark \ref{remark:stop-exploring} for more discussion.
\end{remark}

\begin{prop}\label{prop: terminal string decomposition}
We have the following properties.
\begin{enumerate}
\item For all triples $(s,K,Q) \in \partial_0 \Omega$, $\phi(s,K)=0$.
\item For all triples $(s,K,Q) \in \partial_1 \Omega$, $|\mathrm{supp}(B-K)|+ 8dM(P)\geq \modarea(\ell)$. 
\end{enumerate}
\begin{remark}\label{remark:stop-exploring}
Following Remark \ref{remark:ptl-0-omega}, this proposition says that when we are done exploring, either we realize that there cannot be a surface spanning the loop, and this case contributes zero, or we have seen $\gtrsim \modarea(\ell)$ good or bad plaquettes. Put differently, to reach a terminal string in $\ptl_1 \Omega$ starting from $(\ell, B, \varnothing)$, we must either have taken many deformation steps, or gotten stuck in large bad plaquette clusters. We will pay an exponential cost in the combined total number, which by point (2), is at least given by $\modarea(\ell)$. Thus, the cost will at least be exponential in $\modarea(\ell)$. This is ultimately the reason why we can stop exploring once we have reached the boundary $\ptl \Omega$.
\end{remark}
\end{prop}
\begin{proof}
We begin with the first claim. As previously remarked, if $(s, K)$ contains an isolated edge, then it cannot be balanced, and thus $\phi(s, K) = 0$ by Lemma \ref{lemma:not-balanced-zero}. Similarly, if $s$ is stuck in some minimal union of connected components $\tilde{Q} \sse \mc{P}_\Lambda^+$, and $|\tilde{Q}| < \modarea(s)$, then note by the definition of $\mc{R}(s, K, Q)$ and Proposition \ref{prop:mle2}, we have that
\begin{equs}
\phi(s, K, Q) = \sum_{(s', K', Q') \in \mc{R}(s, K, Q)} c(s', K', Q') \phi(s', K', Q'),
\end{equs}
where here $\phi(s, K, Q) := \phi(s, K)$, and $c(s', K', Q')$ is some combinatorial factor. Thus if $\mc{R}(s, K, Q) = \varnothing$, we obtain that $\phi(s, K) = 0$. We have thus shown the first claim of the proposition. \\

\noindent For the second claim, we split into the following cases. In the following, we assume that $(s, K, Q) \in \ptl_1 \Omega$. \\

\noindent \textit{Case 1:} If $s$ has at least one good edge $\mbf{e}$ and at least one of the merging, splitting, or deformation operations is possible, then since any of these operations increases $|\mathrm{supp}(B-K)\cup Q|$ by at most $1$, we must have that $|\mathrm{supp}(B-K)\cup Q|=\modarea(\ell)$ if $(s,K,Q) \in \partial\Omega$. But $\mathrm{supp}(B-K)\cup Q$ is a connected cluster of $\ell$ (by Lemma \ref{lemma:string-trajectory-cluster}) and $Q$ is solely made up bad plaquettes, thus $|Q| \leq M(P)$ from the definition \ref{eq:M-P} of $M(P)$ and thus $|\mathrm{supp}(B-K)| + M(P) \geq \modarea(\ell)$. \\

\noindent \textit{Case 2:} Suppose that $s$ is stuck in a minimal union of connected components $\tilde{Q}$, and 
\begin{equs}\label{eq:case-2-working-assumption}
|\mathrm{supp}(B-K)\cup Q \cup \tilde{Q}|\leq \modarea(\ell)-8(d-1)M(P).
\end{equs}
This implies that $|Q \cup \tilde{Q}| \leq M(P)$ (by the definition \ref{eq:M-P} of $M(P)$, combined with the fact that $\mrm{supp}(B-K) \cup Q \cup \tilde{Q}$ is a connected cluster of $\ell$ (by Lemma \ref{lemma:string-trajectory-cluster})). Next, by inductively applying Lemma \ref{lemma:area-after-loop-operation-original}, we obtain 
\begin{equs}
\modarea(\ell) &\leq \modarea(s)+|\mathrm{supp}(B-K)|+8(d-1)|Q|.
\end{equs}
By combining the previous considerations, the desired result will now follow upon showing that $\modarea(s) \leq |\tilde{Q}|$. Since $(s, K, Q) \in \ptl_1 \Omega$, by definition, either $\modarea(s) \leq |\tilde{Q}|$ already, or $\mc{R}(s, K, Q) \neq \varnothing$. We proceed show that under the working assumption \eqref{eq:case-2-working-assumption}, the condition $\mc{R}(s, K, Q) \neq \varnothing$ actually implies that $\modarea(s) \leq |\tilde{Q}|$, which will finish the proof of this case. Towards this end, assume that $\mc{R}(s, K, Q) \neq \varnothing$. We argue by contradiction. If $|\tilde{Q}|<\area(s)$, then since $\mc{R}(s,K,Q) \neq \varnothing$, there would exist a nontrivial string trajectory starting at $(s,K,Q)$. Let $(s_1,K_1,Q_1)$ be the next triple in the trajectory, then necessarily $\mathrm{supp}(B-K_1) \sse  \mathrm{supp}(B-K) \cup \partial \tilde{Q}$ and $Q_1=Q \cup \tilde{Q}$, so 
\begin{equs}
|\mathrm{supp}(B-K_1) \cup Q_1| &\leq |\mathrm{supp}(B-K)\cup Q \cup \tilde{Q}|+|\partial \tilde{Q}| \\
&\leq |\mathrm{supp}(B-K)\cup Q \cup \tilde{Q}|+8(d-1)|\tilde{Q}| \\
&\leq \modarea(\ell),
\end{equs}
where in the final inequality we used \eqref{eq:case-2-working-assumption} and $|\tilde{Q}|\leq M(P)$. Thus $(s_1,K_1,Q_1) \in \Omega$, which contradicts the assumption that $(s,K,Q)$ is terminal. This shows that $\modarea(s) \leq |\tilde{Q}|$, finishing the proof of the present case. \\


\noindent
\textit{Case 3:} Suppose that $s$ is stuck in a minimal union of connected components $\tilde{Q} \sse P$ with $|\mathrm{supp}(B-K) \cup \tilde{Q} \cup Q| > \modarea(\ell)-8(d-1) M(P)$. We claim that upon possibly deleting some plaquettes from $\tilde{Q}$, we may obtain a subset $Q' \sse Q \cup \tilde{Q}$ such that $\mrm{supp}(B-K) \cup Q'$ is a connected cluster of $\ell$ with size 
\begin{equs}\label{eq:cluster-size-claim}
\modarea(\ell)-8(d-1) M(P)\leq |\mrm{supp}(B-K) \cup Q'|\leq\modarea(\ell)
\end{equs}
(the main point here is that the upper bound also holds, which is why we may have to delete some plaquettes). We show this claim at the end. Given this claim, 
we now must have $|Q'| \leq M(P)$ (by the definition \eqref{eq:M-P} of $M(P)$), and thus $|\mathrm{supp}(B-K)| \geq \modarea(\ell)-8(d-1)M(P)-M(P)$, and as a consequence, $|\mathrm{supp}(B-K)|+8dM(P) \geq \modarea(\ell)$, completing the proof of this case, modulo the claim.

To show the claim, note by Lemma \ref{lemma:string-trajectory-cluster} that both $\mrm{supp}(B-K) \cup Q$ and $\mrm{supp}(B-K) \cup Q \cup \tilde{Q}$ are connected clusters of $\ell$. Moreover, by the definition of the space $\Omega$ (Definition \ref{def:configuration-space}), we have that $|\mrm{supp}(B-K) \cup Q| \leq \modarea(\ell)$. Our goal is to find a subset $\tilde{Q}' \sse \tilde{Q}$ such that upon setting $Q' = Q \cup \tilde{Q}'$, we have that $\mrm{supp}(B-K) \cup Q'$ is a connected cluster of $\ell$ satisfying \eqref{eq:cluster-size-claim}. Towards this end, first decompose $\tilde{Q} =  \tilde{Q}_1 \cup \cdots \cup \tilde{Q}_k$ into connected components. By our previous observation, $\mrm{supp}(B-K) \cup Q$ and $\mrm{supp}(B-K) \cup Q \cup \tilde{Q}$ are connected, and thus we must have that each component $\tilde{Q}_j$ is connected to $\mrm{supp}(B-K) \cup Q$. Thus for each component $\tilde{Q}_j$, we may take an element $p_j \in \tilde{Q}_j$ which is connected to $\mrm{supp}(B-K) \cup Q$. We may then take a spanning tree of $\tilde{Q}_j$ rooted at $p_j$. We then begin to delete elements of $\tilde{Q}_1$, starting from the leaves of the spanning tree, and working our way up towards the root. If we have completely deleted $\tilde{Q}_1$ and the upper bound \eqref{eq:cluster-size-claim} is still not satisfied, we move onto $\tilde{Q}_2$ and repeat. We stop deleting plaquettes once \eqref{eq:cluster-size-claim} is satisfied, at which point we have arrived at the desired set $\tilde{Q}'$. The point here is that by starting the deletion process from the leaves of the trees, we preserve connectivity after each deletion, so we are guaranteed that the set $\tilde{Q}'$ is such that $\mrm{supp}(B-K) \cup Q \cup \tilde{Q}'$ is a connected cluster of $\ell$.
\\

\noindent \textit{Case 4:} If $s$ is the null loop, then there exists a plaquette count $K'$ supported on $\mrm{supp}(B-K) \cup Q$ such that $(\ell, K')$ is balanced. It follows that
\begin{equs}
|\mrm{supp}(B-K)| + |Q| \geq |\mrm{supp}(K')| \geq \modarea(\ell),
\end{equs}
and thus we may conclude by using that $|Q| \leq M(P)$ (note by definition of $\Omega$ that $|\mrm{supp}(B-K) \cup Q| \leq \modarea(\ell)$, so that $|Q| \leq M(P)$ follows by the definition \eqref{eq:M-P} of $M(P)$ and Lemma \ref{lemma:string-trajectory-cluster}).\\

\noindent This concludes the proof of Proposition \ref{prop: terminal string decomposition}.
\end{proof}

In the following, recall the notion of splitting complexity (Definition \ref{def:splitting-complexity}) which was introduced in Section \ref{section:truncated-model} in the course of the contraction mapping argument. To carry out the contraction mapping argument in the present setting, we will also need to understand how the splitting complexity changes under a revival operation, which is the purpose of the next lemma.

\begin{lemma}[Revival splitting complexity]\label{lemma:splitting-complexity-revival}
Let $s\in \mc{S}_{\Lambda}$ be a non-null stuck string, and let $(s',K',Q') \in \mc{R}(s,K,Q)$. Then $\iota(s')=\iota(s)$.
\end{lemma}

\begin{proof}
Suppose $s'= s \Delta_{\partial \tilde{Q}} J$ for some $\tilde{Q} \sse P$ and plaquette count $J$ with support $\mathrm{supp}(J)=\{p_1,\dots,p_m\}$. Then $s'$ has $|J|=\sum_{i=1}^{m}J(p_i)$ more loops and $|s'|=4|J|+|s|$. Thus $\iota(s')-\iota(s)=(4|J|)/4-|J|=0$.
\end{proof}

Next, we define a norm on the space of functions $f : \Omega \ra \R$.

\begin{definition}[Norm on configuration space]
For parameters $\lambda, \gamma,\rho, > 0$, define the following norm on $\R^{\Omega}$:
\begin{equs}
\|f\|_{\lambda,\gamma, \rho}:=\sup_{(s,K,Q) \in \Omega} \gamma^{|\mathrm{supp}(B-K) \backslash \ptl Q|}\rho^{\sum_{p \in \mc{P}_{\Lambda}^+ \backslash (P \cup \partial Q)} (B-K(p))}\lambda^{\iota(s)}\exp(10^{-2} B|Q|) |f(s,K,Q)|.
\end{equs}
Since $\Omega$ is finite, every function $f \in \R^{\Omega}$ has finite norm.
\end{definition}

\begin{remark}[Conceptual remarks]
We make several conceptual remarks about the norm. First, the parameter $\gamma$ should really be thought of as the inverse of the parameter $\gamma$ from the corresponding norm for the truncated model (Definition \ref{prop:fixed-point-estimate}). In particular, one should think of $\gamma \ll 1$. This is due to the fact that for the initial configuration $(\ell, B , \varnothing)$, the weight $|\mrm{supp}(B-K) \backslash \ptl \varnothing| = |\mrm{supp}(B - B)| = 0$, whereas in the norm of Definition \ref{def:truncated-norm}, the weight starts off at $\area(\ell)$. One should now think of $\gamma^{|\mrm{supp}(B-K) \backslash \ptl Q|}$ as tracking the total ``cost" that we have accumulated during the exploration process. For a technical reason, we want to exclude the plaquettes of $\ptl Q$ from this cost.

Next, the additional weight $\exp(10^{-2} B |Q|)$ is present to provide a bit of extra room when bounding the revival term. This is because in principal, the revival term does not come associated with any small parameters. So the only way to absorb the revival term is to add a weight that increases when $Q$ increases (this is similar in spirit to the weight $\lambda^{\iota(s)}$ which is needed to absorb the splitting term). Thus, we can think of $\exp(10^{-2} B |Q|)$ as tracking the total size of the revivals we have seen. In the end this is sufficient for our purposes, since we have shown in Proposition \ref{prop: reduction to fixed plaquette config} that such a weight can still be absorbed. This is ultimately due to the fact that all plaquettes in $Q$ are bad, and bad plaquettes are extremely rare.
\end{remark}

Next, we define the operator which appears in the modified master loop equation from Section \ref{section:modified-mle}. Recall that $\Omega^\circ = \Omega \backslash \ptl \Omega$.

\begin{definition}[Master loop equation operator]
For $(s,K,Q) \in \Omega^\circ$, if $s$ has at least one good edge, let $\mbf{e}$ be its designated good edge and define
\begin{equs}
(Mf)(s,K,Q):=\mp \sum_{s' \in \mbb{S}_{\pm}(s, \mbf{e})} f(s', K, Q) \mp \frac{1}{N^2} \sum_{s' \in \mbb{M}_{\pm}(s, \mbf{e})} f(s', K, Q) \mp \frac{\beta}{N} \sum_{(s', K') \in \mbb{D}_{\pm}(s, \mbf{e}, K)} f(s', K', Q).
\end{equs}
If $(s,K,Q) \in \Omega^\circ$ is stuck in the minimal union of connected components $Q'$, set 
\begin{equs}
(Mf)(s, K, Q) := \sum_{\substack{ J: P_{\Lambda} \to \N \\ \varnothing \neq \mathrm{supp}(J) \sse \partial Q' \\ J \leq K}} \frac{(N\beta)^J}{J!} f((s,K) \Delta_{\partial Q} J,Q').
\end{equs}
Lastly if $(s,K,Q) \in \partial\Omega$, set $(Mf)(s,K,Q):=\phi(s,K)$.
\end{definition}

We remark that with $\phi(s, K, Q) := \phi(s, K)$, where $\phi(s, K)$ is defined as in \eqref{eq:phi-s-K-original}, we have that $\phi = M \phi$, which folllows from Propositions \ref{prop:mle1} and \ref{prop:mle2}. We are now ready to state and prove the contraction estimate for $M$, which is the main intermediate step towards the proof of Theorem \ref{thm:area-law}.

\begin{prop}[Contraction estimate]\label{prop:fixed-point-estimate}
Let $0<\lambda, \gamma, \rho \leq 1$. Then for $f \in \R^{\Omega}$,
\begin{equs}
\|Mf\|_{\lambda, \gamma, \rho} &\leq \bigg(2dB\lambda+\frac{2dB}{\lambda N^2}+\frac{4d\beta}{\lambda \gamma \rho N} + \exp(-10^{-3}B)\bigg)\|f\|_{\lambda, \gamma, \rho} \\
&\quad \quad +\sup_{(s,K,Q) \in \partial_1 \Omega} \gamma^{|\mathrm{supp}(B-K) \backslash \ptl Q|}\rho^{\sum_{p \in \mc{P}_{\Lambda}^+ \backslash (P \cup \partial Q)} (B-K(p))} \exp(10^{-2}B |Q|) |\phi(s,K)|.
\end{equs}
\end{prop}

\begin{proof}
For notational convenience, given a triple $(s, K, Q)$, define
\begin{equs}
w(s, K, Q) := \gamma^{|\mathrm{supp}(B-K) \backslash \ptl Q|}\rho^{\sum_{p \in \mc{P}_{\Lambda}^+ \backslash (P \cup \partial Q)} (B-K(p))} \lambda^{\iota(s)} \exp(10^{-2}B |Q|).
\end{equs}
For any $(s,K,Q) \in \partial \Omega$, either $(s,K,Q) \in \partial_0 \Omega$ in which case $(Mf)(s,K,Q)=\phi(s,K) = 0$ (by Proposition \ref{prop: terminal string decomposition}), or $(s,K,Q) \in \partial_1 \Omega$, and so
\begin{equs}
\sup_{(s,K,Q) \in \partial \Omega} & w(s, K, Q)  |(Mf)(s,K,Q)| \leq \sup_{(s,K,Q) \in \partial_1 \Omega}  \gamma^{|\mathrm{supp}(B-K) \backslash \ptl Q|}\rho^{\sum_{p \in \mc{P}_{\Lambda}^+ \backslash (P \cup \partial Q)} (B-K(p))} \exp(10^{-2}B |Q|) |\phi(s,K)|,
\end{equs}
where we used that $\lambda < 1$ and $\iota(s) \geq 0$ for any $s \in \mc{S}_{\Lambda}$. To finish, it suffices to show that
\begin{equs}
\sup_{(s,K,Q) \in \Omega^\circ} & w(s, K, Q) |(Mf)(s,K,Q)| \leq \bigg(2dB\lambda+\frac{2dB}{\lambda N^2}+\frac{4d\beta}{\lambda \gamma \rho N} + \exp(-10^{-3}B)\bigg)\|f\|_{\lambda, \gamma, \rho}.
\end{equs}
Towards this end, fix $(s, K,Q) \in \Omega^\circ$. By the definition of $M$, it suffices to prove the following four estimates:
\begin{equs}
w(s, K, Q) \sum_{s' \in \mbb{S}_{\pm}(s, \mbf{e}_s)} |f(s', K, Q)| &\leq 2dB \lambda \|f\|_{\lambda, \gamma, \rho} \label{eq:intermediate-splitting-untruncated}, \\
w(s, K, Q)  \sum_{s' \in \mbb{M}_{\pm}(s, \mbf{e}_s)} |f(s', K)| &\leq \frac{2d B}{\lambda} \|f\|_{\lambda, \gamma, \rho} \label{eq:intermediate-merger-untruncated}, \\ 
w(s, K, Q)  \sum_{(s', K') \in \mbb{D}_{\pm}(s, \mbf{e}_s, K)} |f(s', K')| &\leq  \frac{4d}{\lambda \gamma \rho}  \|f\|_{\lambda, \gamma, \rho}, 
\label{eq:intermediate-deformation-untruncated}\\
w(s, K, Q) \sum_{\substack{(s', K',Q') \in \mc{R}(s, K, Q)\\ s' = s \Delta_{\partial Q'} J}} \frac{(N\beta)^J}{J!}|f(s', K',Q')| &\leq \exp(-10^{-3}B)  \|f\|_{\lambda, \gamma, \rho} 
\label{eq:intermediate-revival-untruncated}.
\end{equs}

\noindent \emph{Proof of \eqref{eq:intermediate-splitting-untruncated}.}
Let $s' \in \mbb{S}_{\pm}(s, \mbf{e}_s)$. By Lemma \ref{lemma:splitting-complexity-after-splitting}, we have that $\iota(s') \leq \iota(s) - 1$, and so $\lambda^{\iota(s)} = \lambda \cdot \lambda^{\iota(s) - 1} \leq \lambda \cdot \lambda^{\iota(s')}$ (since $\lambda \leq 1$). By the definition of splittings, no other part of $w(s, K, Q)$ changes besides $\lambda^{\iota(s)}$, and thus
\begin{equs}
w(s, K, Q) |f(s', K, Q)| \leq\lambda w(s', K, Q)|f(s', K, Q)| \leq \lambda\|f\|_{\lambda,\gamma,\rho}.
\end{equs}
Since $\mbf{e}_s$ is a good edge, there are at most $2dB$ positive or negative splittings at $\mbf{e}_s$ (recall Remark \ref{remark:bound-on-loop-operations-string-trajectory}), and so the estimate \eqref{eq:intermediate-splitting-untruncated} follows. \\

\noindent \emph{Proof of \eqref{eq:intermediate-merger-untruncated}.} 
Let $s' \in \mbb{M}_{\pm}(s, \mbf{e}_s)$. By Lemma \ref{lemma:splitting-complexity-merger} and the fact that $\lambda \leq 1$, we have that $\lambda^{\iota(s)} = \lambda^{-1} \lambda^{\iota(s) + 1} \leq \lambda^{-1} \lambda^{\iota(s')}$. By the definition of splittings, no other part of $w(s, K, Q)$ changes besides $\lambda^{\iota(s)}$, and thus 
\begin{equs}
w(s, K, Q) |f(s', K, Q)| \leq \lambda^{-1 }w(s', K, Q) |f(s', K, Q)|\leq \lambda^{-1}\|f\|_{\lambda,\gamma,\rho}.
\end{equs}
Since $\mbf{e}_s$ is a good edge, there are at most $2dB$ positive or negative mergers (recall Remark \ref{remark:bound-on-loop-operations-string-trajectory}), and so the estimate \eqref{eq:intermediate-merger-untruncated} follows. \\

\noindent \emph{Proof of \eqref{eq:intermediate-deformation-untruncated}.}
Let $(s', K') \in \mbb{D}_{\pm}(s, \mbf{e}_s, K)$. By Lemma \ref{lemma:splitting-complexity-merger}, $\iota(s') \leq \iota(s) + 1$, and so $\lambda^{\iota(s)} \leq \lambda^{-1} \lambda^{\iota(s')}$ (as in the merger case). Since $K'$ is obtained from $K$ by subtracting 1 from a single plaquette value $K(p)$, and $Q' = Q$, we also have that $|\mathrm{supp}(B-K') \backslash \ptl Q'| \leq |\mathrm{supp}(B-K) \backslash \ptl Q|+1$, and since $\gamma \leq 1$, this implies that
\begin{equs}
\gamma^{|\mrm{supp}(B-K) \backslash \ptl Q|} \leq \gamma^{-1} \gamma^{|\mrm{supp}(B-K') \backslash \ptl Q'|}.
\end{equs}
Similarly, we also have that
\begin{equs}
\sum_{p \in \mc{P}_{\Lambda}^+ \backslash (P\cup \partial Q)} (B-K'(p)) \leq \sum_{p \in \mc{P}_{\Lambda}^+ \backslash (P\cup \partial Q)} (B-K(p)) +1,
\end{equs}
and thus since $\rho \leq 1$, we obtain
\begin{equs}
\rho^{\sum_{p \in \mc{P}_{\Lambda}^+ \backslash (P\cup \partial Q)} (B-K(p))} \leq \rho^{-1} \rho^{\sum_{p \in \mc{P}_{\Lambda}^+ \backslash (P\cup \partial Q)} (B-K'(p))}.
\end{equs}
Combining these considerations, we obtain
\begin{equs}
w(s, K, Q) |f(s', K', Q')|&\leq \lambda^{-1} \gamma^{-1}  \rho^{-1} w(s', K', Q) |f(s', K', Q')| \\
&\leq \lambda^{-1} \gamma^{-1} \rho^{-1}\|f\|_{\lambda,\gamma,\rho}.
\end{equs}
Since there are at most $4d$ positive or negative deformations (by Lemma \ref{lemma:bound-on-loop-operations}), the estimate \eqref{eq:intermediate-deformation-untruncated} follows. \\

\noindent \emph{Proof of \eqref{eq:intermediate-revival-untruncated}.}
Let $(s', K',Q') \in \mc{R}(s,K,Q)$. By Lemma \ref{lemma:splitting-complexity-revival}, we have that $\iota(s') =\iota(s)$. Also, by the definition of $\mc{R}(s, K, Q)$, $K'$ is obtained from $K$ by zeroing out the plaquettes of $\ptl(Q' \backslash Q) \sse \ptl Q'$ (recall Remark \ref{remark:edge-total-monotone}), and thus we have that $\mrm{supp}(B-K') \backslash \ptl Q' \sse \mrm{supp}(B-K) \backslash \ptl Q$.
Since $\gamma \leq 1$, we obtain
\begin{equs}
\gamma^{|\mathrm{supp}(B-K)\backslash \ptl Q|} \leq \gamma^{|\mathrm{supp}(B-K')\backslash \ptl Q'|}.
\end{equs}
Similarly, by the definition of $\mc{R}(s, K, Q)$, we have that $K'=K$ on $ \mc{P}_{\Lambda}^+ \backslash (P \cup \partial Q')$ (again since $K'$ is obtained from $K$ by only modifying its values on $\ptl (Q' \backslash Q) \sse \ptl Q'$), we have that (using also that $\ptl Q \sse \ptl Q'$)
\begin{equs}
\sum_{p \in \mc{P}_{\Lambda}^+ \backslash (P \cup \partial Q)} (B-K(p)) \geq \sum_{p \in \mc{P}_{\Lambda}^+ \backslash (P \cup \partial Q')} (B-K'(p)).
\end{equs}
Since $\rho \leq 1$, we obtain
\begin{equs}
\rho^{\sum_{p \in \mc{P}_{\Lambda}^+ \backslash (P \cup \partial Q)} (B-K(p))} \leq \rho^{\sum_{p \in \mc{P}_{\Lambda}^+ \backslash (P \cup \partial Q')} (B-K'(p))}.
\end{equs}
Combining these considerations, we obtain
\begin{equs}
w(s, K, Q) |f(s', K',Q')| &\leq \exp(-10^{-2}B(|Q'|-|Q|)) w(s', K', Q') |f(s', K',Q')|\\
&\leq \exp(-10^{-2}B(|Q'|-|Q|))\|f\|_{\lambda,\gamma,\rho},
\end{equs}
and as a result
\begin{equs}
w(s, K, Q) &\sum_{\substack{(s', K',Q') \in \mc{R}(s, K, Q)\\ s' = s \Delta_{\partial Q'} J}} \frac{(N\beta)^J}{J!}|f(s', K',Q')|\\
&\leq \exp(-10^{-2}B(|Q'|-|Q|))\|f\|_{\lambda,\gamma,\rho} \sum_{\substack{(s', K',Q') \in \mc{R}(s, K, Q)\\ s' = s \Delta_{\partial Q'} J}} \frac{(N\beta)^J}{J!}\\
&\leq \exp(-10^{-2}B(|Q'|-|Q|))\|f\|_{\lambda,\gamma,\rho} \exp(16 N\beta (d-1) (|Q'|-|Q|))\\
&\leq \exp(-10^{-3}B(|Q'|-|Q|))\|f\|_{\lambda,\gamma,\rho}\\
&\leq \exp(-10^{-3}B) \|f\|_{\lambda,\gamma,\rho}.
\end{equs}
Here, we used that due to our choice of parameters, we have that $16N\beta (d-1) \ll B$ and $\exp(-10^{-3} B) < 1$. Also, in the last step, we used the fact that $|Q'| - |Q| \geq 1$.
This concludes the proof of \eqref{eq:intermediate-revival-untruncated} and thus also the proof of Proposition \ref{prop:fixed-point-estimate}.
\end{proof}

Proposition \ref{prop:fixed-point-estimate} directly leads to the following corollary under our assumptions on the parameters $N, \beta, B$.

\begin{cor}\label{cor:Fixed-point-with-well-chosen-parameters}
Under the choice of parameters \eqref{eq:N-large-assumption}-\eqref{eq:B-choice}, let $\lambda=N^{-1}$, $\gamma=10^3 d \beta \leq 1$, $\rho=e^{-1}$.  We have that
\begin{equs}
\|Mf\|_{\lambda, \gamma, \rho}\leq \frac{1}{2}\|f\|_{\lambda, \gamma, \rho} +\sup_{(s,K,Q) \in \partial_1 Q} \gamma^{|\mathrm{supp}(B-K) \backslash \ptl Q|}\rho^{\sum_{p \in \mc{P}_{\Lambda}^+ \backslash (P \cup \partial Q)}(B-K(p))}\exp(10^{-2}B|Q|)|\phi(s,K,Q)|.
\end{equs}
\end{cor}
\begin{proof}
This follows by Proposition \ref{prop:fixed-point-estimate}, combined with the fact that
\begin{equs}
2dB\lambda+\frac{2dB}{\lambda N^2}+\frac{4d\beta}{\lambda \gamma \rho N} + \exp(-10^{-3}B) \leq \frac{1}{2}.
\end{equs}
To see this, use that $\lambda = N^{-1}$, $\rho = e^{-1}$, and recall that $dB \ll N$ (which handles the first two terms), 
\begin{equs}
\frac{d \beta}{\gamma} = 10^{-3} \ll 1,
\end{equs}
which handles the third term, and $\exp(-10^{-3} B) \ll 1$, which handles the last term.
\end{proof}

Before we prove Theorem \ref{thm:area-law}, we estimate the terminal string integrals appearing in Corollary \ref{cor:Fixed-point-with-well-chosen-parameters}.

\begin{lemma}[Bound on terminal string integrals]\label{lemma:terminal-string-bound}
For all $(s,K,Q) \in \partial_1 \Omega$, we have that
\begin{equs}
\exp(10^{-2} B |Q|)|\phi(s,K)| \leq \exp((B/50) |Q|) e^{\sum_{p\in \mc{P}_{\Lambda}^+ \backslash (P\cup \partial Q)}(B-K(p))} Z_{\Lambda,P},
\end{equs}
where 
\begin{equs}
Z_{\Lambda,P}:=\int dU \prod_{p \in \mc{P}_{\Lambda}^+ \backslash P} \exp_B(2\beta \mathrm{Re}\Tr(U_p))\prod_{p \in  P} \tail_B(2\beta \mathrm{Re}\Tr(U_p))).
\end{equs}
\end{lemma}
\begin{proof}
We trivially have the pointwise bound $|W_s(U)| \leq 1$ for all $U : E_\Lambda^+ \ra \unitary(N)$. Combining this with \eqref{eq:Truncation level comparison}, we have that
\begin{equs}
|\phi(s,K)|&=\bigg|\int dU W_s(U) \prod_{p \in \mc{P}^+_{\Lambda}\backslash P}\exp_{K(p)}(2\beta \mathrm{Re} \mathrm{Tr}(U_p))\prod_{p \in P} \tail_B(2\beta \mathrm{Re} \mathrm{Tr}(U_p))) \bigg|\\
&\leq \int dU \prod_{p \in \mc{P}^+_{\Lambda}\backslash P} \exp_{K(p)}(2\beta \mathrm{Re} \mathrm{Tr}(U_p)) \prod_{p \in P} \tail_B(2\beta \mathrm{Re} \mathrm{Tr}(U_p))) \\
&\leq \int dU \prod_{p \in \partial Q} e^{2N\beta} \prod_{p \in \mc{P}^+_{\Lambda}\backslash (P \cup \ptl Q)}e^{B-K(p)}\exp_{B}(2\beta \mathrm{Re} \mathrm{Tr}(U_p))\prod_{p \in P} \tail_B(2\beta \mathrm{Re} \mathrm{Tr}(U_p))).
\end{equs}
Here, we applied the bound $\exp_{K(p)} \leq e^{B-K(p)} \exp_B$ for the terms $p \in \mc{P}_\Lambda^+ \backslash (P \cup \ptl Q)$ (recall \eqref{eq:Truncation level comparison}), and the bound $|\exp_{K(p)}(x)| \leq \exp(2N\beta)$ for all $|x| \leq 2N\beta$ for the terms $p \in \ptl Q$. Next, we use that (by \eqref{eq:Truncation level comparison}) for all $|x| \leq 2N\beta$,
\begin{equs}
\exp_B(x) = \exp(x) - \tail_B(x) \geq e^{-2N\beta} - e^{-B/10} \exp_B(x),
\end{equs}
so that (since $B \gg 1$ so $e^{-B/10} \ll 1$)
\begin{equs}
\exp_B(x) \geq \frac{1}{2} e^{-2N\beta},
\end{equs}
and thus we obtain the further upper bound
\begin{equs}
|\phi(s, K)| \leq 2^{|\ptl Q|} e^{4N\beta |\ptl Q|} e^{\sum_{p \in \mc{P}_\Lambda^+ \backslash (P \cup \ptl Q)} (B - K(p))} Z_{\Lambda, P}.
\end{equs}
To finish, we use that $|\ptl Q| \leq Cd|Q|$ and $d N \beta \ll B$, $d \ll B$, which holds due to our choice of parameters.
\end{proof}

We are finally ready to prove Theorem \ref{thm:area-law}.

\begin{proof}[Proof of Theorem \ref{thm:area-law}]
Let $f(s, K, Q) := \phi(s, K)$. By Propositions \ref{prop:mle1} and \ref{prop:mle2}, we have that $Mf = f$, and thus by Corollary \ref{cor:Fixed-point-with-well-chosen-parameters}, upon setting $\lambda = N^{-1}$, $\gamma = 10^3 d \beta$, $\rho = e^{-1}$, we obtain that
\begin{equs}
\|f\|_{\lambda, \gamma, \rho}&\leq 2\sup_{(s,K,Q) \in \partial_1 \Omega} \gamma^{|\mathrm{supp}(B-K) \backslash \ptl Q|}\rho^{\sum_{p \in P_{\Lambda}^+ \backslash (P \cup \partial Q)}(B-K(p))}\exp(10^{-2}B |Q|)|\phi(s,K,Q)|\\
&\leq 2 \sup_{(s, K, Q) \in \ptl_1 \Omega} \gamma^{|\mrm{supp}(B-K) \backslash \ptl Q|} \exp((B/50)|Q|) Z_{\Lambda, P}, \\
&\leq 2 \sup_{(s, K, Q) \in \ptl_1\Omega} \gamma^{(\modarea(\ell) - 8dM(P) - |\ptl Q|)_+} 
\exp((B/50)  |Q|)Z_{\Lambda,P}.
\end{equs}
where we used Lemma \ref{lemma:terminal-string-bound} in the second inequality, and Proposition \ref{prop: terminal string decomposition} in the third inequality (combined with $\gamma \leq 1$). Next, using that $|Q| \leq M(P)$ (since $\mrm{supp}(B-K) \cup Q$ is a connected cluster by Lemma \ref{lemma:string-trajectory-cluster}, and moreover its size is at most $\modarea(\ell) = \area(\ell)$ by the definition of $\Omega$), which also implies that $|\ptl Q| \leq 8d M(P)$, we further obtain that
\begin{equs}
\|f\|_{\lambda, \gamma, \rho} \leq 2 \gamma^{(\modarea(\ell) - 20 d M(P))_+} \exp((B/50)M(P)) Z_{\Lambda, P}.
\end{equs}
In summary, we obtain (using that $\area(\ell) = \modarea(\ell)$ for rectangular loops, $\iota(\ell) = |\ell|/4 - 1$, and that $\lambda = N^{-1}$)
\begin{equs}
\phi(\ell,B) \leq N^{|\ell|/4-1} \|f\|_{\lambda,\rho,\gamma} \leq 2 N^{|\ell|/4-1} \times (10^3 d\beta)^{(\area(\ell) - 20dM(P))_+} \exp((B/40)|M(P)|) Z_{\Lambda,P}.
\end{equs}
Plugging this into Proposition \ref{prop: reduction to fixed plaquette config}, we obtain (using that $|\ell| \leq \area(\ell)$, so that $e^{2d|\ell|} \leq e^{2d|\area(\ell)|}$)
\begin{equs}
|\langle W_\ell \rangle_{\Lambda, \beta, N}| &\leq 2|\ell| e^{2d|\ell|} N^{|\ell|/4-1} \times 50^{d\cdot \area(\ell)} 2^{\area(\ell)} \max(10^3 d \beta, \exp(-10^{-3} d^{-1} B))^{\area(\ell)} \\
&\leq 2 |\ell| N^{|\ell|/4-1} \Big(2 \times 10^{3d} \max(10^3 d \beta, \exp(-10^{-3} d^{-1} B))\Big)^{\area(\ell)}.
\end{equs}
To finish, we use that $B \geq \frac{1}{2d} 10^{-3} N$ by \eqref{eq:B-choice}, and thus
\[ 
\exp(-10^{-3} d^{-1} B) \leq \exp(- 10^{-7} d^{-2} N). \qedhere
\]
\end{proof}

\bibliographystyle{alpha}
\bibliography{references}

\end{document}